\documentclass[11pt,leqno]{article}

\usepackage{amsmath}
\usepackage{amssymb}
\usepackage{amscd}
\usepackage{graphicx}
\input{epsf}

\newtheorem{Theorem}{Theorem}[section]
\newtheorem{Corollary}[Theorem]{Corollary}
\newtheorem{Lemma}[Theorem]{Lemma}

\newtheorem{Example}[Theorem]{Example}
\newtheorem{maintheorem}{Theorem}

\,
\,
\,

\def\qed{\hfill $\bullet$}
\def\proof{\noindent {\em Proof:}\ }

\def\set-up{\noindent {\em Set-Up:}\ }

\newcommand{\pres}[2]{\langle {#1}\ |\ {#2} \rangle}

\newcommand{\sgp}[1]{\langle {#1}\rangle}

\def\Z{\mathbb Z}

\def\ra{\rightarrow}

\def\<{\langle}
\def\>{\rangle}

\def\'{\,^\prime}
\def\"{\,^{\prime\prime}}

\def\iff{\Leftrightarrow}

\def\ker{\mathrm{ker}\,}

\begin{document}

\title{Cyclically Presented Groups With Length Four Positive Relators}
\author{William A.~Bogley and Forrest W.~Parker}
\date{\today}

\maketitle

\begin{abstract}
We classify cyclically presented groups of the form $G = G_n(x_0x_jx_kx_l)$ for finiteness and, modulo two unresolved cases, we classify asphericity for the underlying presentations. We relate finiteness and asphericity to the dynamics of the shift action by the cyclic group of order $n$ on the nonidentity elements of $G$ and show that the fixed point subgroup of the shift is always finite.
\end{abstract}

\section{Cyclically Presented Groups and Shift Dynamics}\label{Section:Introduction}

One purpose of this article is to provide classifications of finiteness and asphericity (modulo two unresolved cases) for cyclic group presentations where the relators are positive words of length four. More generally, we relate finiteness and asphericity to the dynamics of the shift action that arises from the cyclic symmetry in these presentations, thus extending results from \cite{BShift}.

Let $n$ be a positive integer and let $F$ be the free group of rank $n$ with basis $x_0, \ldots, x_{n-1}$. The \textbf{shift automorphism} is given by $\theta_F(x_i) = x_{i+1}$ where subscripts are taken modulo $n$. A word $w \in F$ determines the \textbf{cyclic (group) presentation}
$$
\mathcal{P}_n(w) = \pres{x_0, \ldots, x_{n-1}}{w,\theta_F(w), \ldots, \theta_F^{n-1}(w)}
$$
which in turn defines the \textbf{cyclically presented group} $G = G_n(w)$. The shift determines an automorphism $\theta_G \in \mathrm{Aut}(G)$ that satisfies $\theta_G^n = 1$. The order of the shift is therefore a divisor of $n$ and the shift may or may not be an inner autormorphism of $G$. 

The following two classification problems have received attention in the literature.  See \cite{CRS05,EdjvetWilliams,TRevisit,GW09}, for example.

\begin{description}
\item[Finiteness:] For which $n$ and $w$ is the cyclically presented group $G_n(w)$ finite?
\item[Asphericity:] For which $n$ and $w$ is the cyclic presentation $\mathcal{P}_n(w)$ (combinatorially) aspherical?
\end{description}
See Section \ref{Section:Asphericity} for a discussion of combinatorial asphericity. A group presentation is \textbf{aspherical} if its two-dimensional cellular model $K$ is aspherical, or equivalently its second homotopy module $\pi_2K$ is trivial. Combinatorial asphericity is a more general concept that allows for the presence of proper power and freely redundant relators, in which case $\pi_2K \neq 0$. See \cite{CCH} for a comprehensive treatment of asphericity concepts for group presentations. Asphericity and finiteness are related in that any finite subgroup of a group defined by a combinatorially aspherical presentation is finite cyclic \cite{Hueb}. This generalizes the well-known fact that if a two-complex $K$ satisfies $\pi_2K = 0$, then the fundamental group $\pi_1K$ is torsion-free.

In \cite{CRS05,EdjvetWilliams}, the finiteness and asphericity classification problems were completely solved for cyclic presentations $\mathcal{P}_n(w)$ where the defining relator is a positive word of length three. In this article we deal with the length four case. Thus we consider presentations of the form 
\begin{equation}\label{Equation:P}
\mathcal{P}_n(x_0x_jx_kx_l) = \pres{x_0, \ldots x_{n-1}}{x_ix_{i+j}x_{i+k}x_{i+l}, 0 \leq i < n}
\end{equation}
where the integer parameters $j,k,l$ are to be considered modulo the positive integer $n$. (Note that the presentation $\mathcal{P}_n(w)$ is unchanged if we replace $w$ by any of its shifts.) Modeling an approach developed in \cite{EdjvetWilliams}, our classifications are framed in terms of conditions on the parameters $(n,j,k,l)$ starting with (A), (B), (C) as follows. All congruences are considered modulo $n$. 

\begin{table}[ht]
\begin{center}
\begin{tabular}{ll}
(A) & $2k \equiv 0$ or $2j \equiv 2l$\\
(B) & $k \equiv 2j$ or $k \equiv 2l$ or $j+l \equiv 2k$ or $j+l \equiv 0$\\
(C) & $l \equiv j + k$ or $l \equiv j-k$\\
\end{tabular}
\end{center}
\caption{Conditions (A), (B), (C)}\label{Table:ABC}
\end{table}
%
%\begin{enumerate}
%\item[\bf (A)]  $2k \equiv 0$ or $2j \equiv 2l$.
%\item[\bf (B)]  $k \equiv 2j$ or $k \equiv 2l$ or $j+l \equiv 2k$ or $j+l \equiv 0$.
%\item[\bf (C)] $l \equiv j + k$ or $l \equiv j-k$.
%\end{enumerate}
\noindent Table \ref{Table:Results} provides an overview of our main results regarding asphericity and finiteness in terms of truth values for the conditions (A), (B), (C).

\begin{table}[ht]
\begin{center}
\begin{tabular}{|c||c|c|c|}\hline
\begin{tabular}{rcl} (A)&(B)&(C) \\  \end{tabular} & Aspherical [\ref{Theorem:AsphericityClassification}] & Group [\ref{Theorem:FinitenessClassification}] &  $W \mbox{[\ref{Lemma:Cdecomp}-\ref{Lemma:B}]}$ \\ \hline\hline
\begin{tabular}{ccc}T& T & T \\ T & F &T\\ T&F&F\\F&F&F\\ \end{tabular} &\begin{tabular}{c}Yes\\C(4)-T(4)\\ $\mbox{[\ref{Lemma:SmCanc}, \ref{Lemma:Decomposable}]}$ \\ \end{tabular} & Finite $\iff n = 1 \iff \Z_4$ [\ref{Lemma:FiniteAsph}(b)]  & \begin{tabular}{c}$z^4$\\ $(z^2\alpha)^2$\\ $\prod_{i=1}^4 x\alpha_i$\\ $\prod_{i=1}^4 x\alpha_i$\\ \end{tabular}\\ \hline
\begin{tabular}{ccc}F& T & T \\ F & F &T\\ \end{tabular} & No [\ref{Lemma:Decomposable}] & \begin{tabular}{c}Finite $\iff \Z_4$ [\ref{Theorem:FiniteCnotA}] \\ See [\ref{Theorem:FinitenessClassification}(b)],[\ref{Theorem:FiniteCnotA}(iv)]\\ \end{tabular} & \begin{tabular}{c}$z^4\alpha$ \\ $z^2\alpha z^2 \beta$\\ \end{tabular}\\ \hline
\begin{tabular}{ccc} T&T&F \\  \end{tabular} &No [\ref{Theorem:MZ}] & \begin{tabular}{c} Finite $\iff \gamma = 1$ \\ Finite $\Rightarrow$ Solvable\\ \end{tabular} [\ref{Theorem:MZ}] & $u^3\alpha u \alpha^{\pm 1}$\\ \hline
\begin{tabular}{ccc} F&T&F \\  \end{tabular}&\begin{tabular}{c}(I$^\ast$): No [\ref{Theorem:INot}]\\ (U$^\ast$): Unresolved\\ Else: Yes [\ref{Theorem:Asphericity}]\\ \end{tabular}& Finite $\iff$ (I5), (I6$'$), or (I6$''$) [\ref{Theorem:FiniteBnotAC}]& $u^3\alpha u \beta$\\ \hline
%&& &(I) $\Rightarrow$ Not CA & (I6$'$) or (I6$''$) & & \\ \hline
\end{tabular}
\end{center}
\caption{Asphericity and Finiteness for $\mathcal{P}_n(x_0x_jx_kx_l)$}\label{Table:Results}
\end{table}

%Finite $\stackrel{[\ref{Theorem:FiniteBnotAC}]}{\iff}$ \begin{tabular}{l} $\gcd(n,j,k,l) = 1$ and \\ (I5), (I6$'$), or (I6$''$)\\ \end{tabular}

%\begin{tabular}{rl} $\Z_4$ & if $n$ odd\\ V1D & if $n$ even\\ \end{tabular}
%
%\begin{tabular}{c} Not CA $\Rightarrow$ (I) or (U)\\ (I) $\Rightarrow$ Not CA\\ \end{tabular}
%
%\begin{tabular}{l} (I5), (I6$'$), \\ or (I6$''$)\\ \end{tabular}

%\begin{tabular}{c}T\\F\\ \end{tabular} & \begin{tabular}{c}F\\F\\ \end{tabular} & \begin{tabular}{c}F\\F\\ \end{tabular}

%\begin{tabular}{c}$u^4$\\$(u^2\alpha)^2$\\ \end{tabular}
%
%\begin{tabular}{ccc}\ T&\ T&\ T \\ T&F&T\\ \end{tabular}
%\begin{tabular}{ccc} T&F&F \\ F&F&F\\ \end{tabular}
%\begin{tabular}{ccc} F&T&T \\  \end{tabular}
%\begin{tabular}{ccc} F&F&T \\  \end{tabular}
%\begin{tabular}{ccc} T&T&F \\  \end{tabular}
%\begin{tabular}{ccc} F&T&F \\  \end{tabular}
%\begin{tabular}{ccc} (A)&(B)&(C) \\  \end{tabular}
%
%
%\begin{tabular}{ccc} $\Z_4$ ($n=1$) or VF \\ \hline $1 \leq \mathrm{vgd} \leq 2$\\ \end{tabular}

%\begin{tabular}{rl} $\Z_4$ & if $n=1$\\ $1 \leq \mathrm{vgd} \leq 2$ & if $n > 1$\\ \end{tabular}
%
%\begin{tabular}{ccc} T&T&T \\ T&F&T\\ \end{tabular}
\noindent We now discuss the specialized notation in Table \ref{Table:Results}. Bracketed entries [-] refer to supporting results in this article. The cyclic group of order $n$ acts on a cyclically presented group $G = G_n(w)$ by automorphisms and so there is the \textbf{shift extension}
$$
E_n(w) = G_n(w) \rtimes_{\theta_G} \Z_n.
$$
This group admits a two-generator two-relator presentation $E_n(w) \cong \pres{a,x}{a^n, W}$ where $W$ is obtained from $w$ via the substitutions $x_i = a^ixa^{-i}$, so the shift $\theta_G \in \mathrm{Aut}(G)$ arises via conjugation by $a$ in $E_n(w)$. When $w = x_0x_jx_kx_l$, the relator $W$ takes the form
\begin{equation}\label{Equation:W}
W = xa^jxa^{k-j}xa^{l-k}xa^{-l}.
\end{equation}
The conditions (A), (B), (C) enable simplifications of the relation (\ref{Equation:W}) as shown in the right hand column of Table \ref{Table:Results}. For example, if $k \equiv 2j$ mod $n$, as in (B), then the substitution $u = xa^j$ transforms $W$ to $u^3a^{l-3j}ua^{-l-j}$. Simplifications of this sort are developed in Section \ref{Section:ABC}. Section \ref{Section:Asphericity} discusses the concept of combinatorial asphericity and its relation to the \textbf{small cancellation} property C(4)-T(4), as well as asphericity for the \textbf{relative presentation} $\pres{\Z_n,x}{W}$ of the shift extension. Section \ref{Section:DU} discusses isomorphisms among groups of the form $G_n(x_0x_jx_kx_l)$ that arise from an action of a finite group $\Gamma_n$ on the set of length four positive words in the generators $x_0, \ldots, x_{n-1}$. The \textbf{secondary divisor} $\gamma = \gcd(n,k-2j,j-2k+l,k-2l,j+l)$ is introduced in Section \ref{Section:Divisor}. 

The various combinations of conditions represented in the lefthand column of Table \ref{Table:Results} are analyzed individually in Sections \ref{Section:SmallCancellation}-\ref{Section:FTF}. As in \cite{BShift,CRS05,EdjvetWilliams}, we depend on previous work concerning asphericity of relative presentations, in this case \cite{BBP}. As is typical in the study of relative asphericity, most of the presentations encountered in \cite{BBP} are aspherical, with the nonaspherical cases falling into well-defined infinite families (see Theorem \ref{Theorem:BBPUpdate}.1), along with a handful of isolated nonaspherical cases. In addition, there are two cases treated in \cite{BBP}  where the asphericity status remains unresolved even when we incorporate recent progress due to Aldwaik and Edjvet \cite{AE14} and Williams \cite{GW15}. See Section \ref{Section:FTF} for these updates.

The bottom row of Table \ref{Table:Results} is the most complex, where we encounter eight \textbf{isolated} conditions and the two \textbf{unresolved} conditions, which are collectively referred to as (I$^\ast$) and (U$^\ast$) respectively. As we now describe, these conditions are defined in terms of \textbf{exemplars} that are listed in Tables \ref{Table:I} and \ref{Table:U}. 
\begin{table}[ht]
\begin{center}
\begin{tabular}{|c|c|c|}\hline
Condition & Exemplar & Shift Extension\\ \hline \hline
(I5) & $\mathcal{P}_5(x_0x_3x_1x_1)$ & $\pres{a,x}{a^5, xa^3xa^{-2}xxa^{-1}} \stackrel{u=xa^3}{\cong} \pres{a,u}{a^5, u^3a^2ua}$ \\ \hline
(I6$'$) & $\mathcal{P}_6(x_0x_4x_2x_3)$ & $\pres{a,x}{a^6, xa^4xa^{-2}xaxa^{-3}} \stackrel{u=xa^4}{\cong} \pres{a,u}{a^6, u^3a^{-3}ua^{-1}}$\\ \hline
(I6$''$) & $\mathcal{P}_6(x_0x_0x_1x_2)$ & $\pres{a,x}{a^6, xxaxaxa^{-2}} \stackrel{u=xa}{\cong} \pres{a,u}{a^6, u^3a^{-3}ua^{-1}}$ \\\hline
(I10) & $\mathcal{P}_{10}(x_0x_3x_6x_1)$ & $\pres{a,x}{a^{10}, xa^3xa^3xa^{-5}xa^{-1}} \stackrel{u=xa^3}{\cong} \pres{a,u}{a^{10}, u^3a^{-8}ua^{-4}}$ \\\hline
(I12) & $\mathcal{P}_{12}(x_0x_1x_2x_9)$ & $\pres{a,x}{a^{12}, xaxaxa^7xa^{-9}} \stackrel{u=xa}{\cong} \pres{a,u}{a^{12}, u^3a^{6}ua^{2}}$ \\\hline
(I16) & $\mathcal{P}_{16}(x_0x_3x_6x_1)$ & $\pres{a,x}{a^{16}, xa^3xa^3xa^{-5}xa^{-1}} \stackrel{u=xa^3}{\cong} \pres{a,u}{a^{16}, u^3a^{-8}ua^{-4}}$ \\\hline
(I20) & $\mathcal{P}_{20}(x_0x_3x_6x_1)$ & $\pres{a,x}{a^{20}, xa^3xa^3xa^{-5}xa^{-1}} \stackrel{u=xa^3}{\cong} \pres{a,u}{a^{20}, u^3a^{-8}ua^{-4}}$ \\\hline
(I24) & $\mathcal{P}_{24}(x_0x_1x_2x_{15})$ & $\pres{a,x}{a^{24}, xaxaxa^{13}xa^{-15}} \stackrel{u=xa}{\cong} \pres{a,u}{a^{24}, u^3a^{12}ua^{8}}$ \\\hline
\end{tabular}
\end{center}
\caption{Exemplars for the Eight Isolated Conditions (I5)-(I24)}\label{Table:I}
\end{table}

\begin{table}[ht]
\begin{center}
\begin{tabular}{|c|c|c|}\hline
Condition & Exemplar & Shift Extension\\ \hline \hline
(U24$'$) & $\mathcal{P}_{24}(x_0x_3x_6x_1)$ & $\pres{a,x}{a^{24}, xa^3xa^3xa^{-5}xa^{-1}} \stackrel{u=xa^3}{\cong} \pres{a,u}{a^{24},u^3a^{-8}ua^{-4}}$\\ \hline
(U24$''$) & $\mathcal{P}_{24}(x_0x_1x_2x_{19})$ & $\pres{a,x}{a^{24}, xaxaxa^{17}xa^{-19}} \stackrel{u=xa}{\cong} \pres{a,u}{a^{24}, u^3a^{16}ua^{4}}$ \\\hline\end{tabular}
\end{center}
\caption{Exemplars for the Two Unresolved Conditions (U24$'$) and (U24$''$)}\label{Table:U}
\end{table}
%
%Each of these ten conditions involves a fixed value of $n$ and cyclic presentations $\mathcal{P}_n(w)$ where $w$ is taken from a single orbit of positive length four words under the action of the finite group $\Gamma_n$. Table \ref{Table:I} displays a single orbit representative, or exemplar, for each of the eight isolated conditions. Table \ref{Table:U} displays exemplars for the unresolved types. 

\paragraph{Types (I) and (U)} As discussed in Section \ref{Section:DU}, the orbit of a positive word $x_ix_jx_kx_l$ under the action of the finite group $\Gamma_n$ consists of all shifts of cyclic permutations of words of the form $x_{iu}x_{ju}x_{ku}x_{lu}$ or $x_{iu}x_{lu}x_{ku}x_{ju}$ where $u \in \Z_n^\ast$ is a multiplicative unit modulo $n$. (The order of the group $\Gamma_n$ is $8n\phi(n)$ where $\phi$ is the Euler totient function.)
%\begin{itemize}
%\item all shifts of cyclic permutations of words of the form $x_{iu}x_{ju}x_{ku}x_{lu}$ or $x_{iu}x_{lu}x_{ku}x_{ju}$ where $u \in \Z_n^\ast$ is a multiplicative unit modulo $n$.
%\end{itemize}
Applying this to the exemplar listed under (I5) in Table \ref{Table:I}, we say that a cyclic presentation $\mathcal{P}_n(w)$ is of \textbf{type (I5)} if and only if $n=5$ and $w$ is a shift of a cyclic permutation of a word of the form $x_0x_{3u}x_ux_u$ or $x_0x_{u}x_ux_{3u}$ where $u \in \Z_5^\ast$. Thus the presentations of type (I5) are precisely those of the form $\mathcal{P}_5(w)$ where $w$ is a length four positive word appearing in the orbit of $x_0x_3x_1x_1$ under the action of the finite group $\Gamma_5$. In Section \ref{Section:DU} it is shown that each such orbit gives rise to a single group up to isomorphism, asphericity status, and shift dynamics. In the same way, the remaining seven exemplars listed in Table \ref{Table:I} determine  corresponding isolated types (I6$^\prime$), (I6$^{\prime\prime}$), (I10), (I12), (I16), (I20), and (I24). Any cyclic presentation that adheres to any one of these eight isolated types is said to be of \textbf{type (I)}. In the same way, the exemplars listed in Table \ref{Table:U} determine the presentations of type (U24$^\prime$) and (U24$^{\prime\prime}$) in terms of the action of the finite group $\Gamma_{24}$. Any presentation of type (U24$^\prime$) or (U24$^{\prime\prime}$) is of \textbf{type (U)}.

\paragraph{Types (I$^\ast$) and (U$^\ast$)} It is routine to verify that every presentation $\mathcal{P}_n(x_0x_jx_kx_l)$ of type (I) or (U) satisfies $\gcd(n,j,k,l) = 1$. For arbitrary $(n,j,k,l)$, we say that $\mathcal{P}_n(x_0x_jx_kx_l)$ is of \textbf{type (I$^\ast$)} (resp.~\textbf{(U$^\ast$)}) if upon setting $c = \gcd(n,j,k,l)$, the presentation $\mathcal{P}_{n/c}(x_0x_{j/c}x_{k/c}x_{l/c})$ is of type (I) (resp.~(U)). The presentations of type (I$^\ast$) are all non-aspherical (see Lemma \ref{Lemma:FreeProduct} and Theorem \ref{Theorem:INot}), but is unknown whether the presentations of type (U$^\ast$) are aspherical. Resolution of this ambiguity reduces to consideration of the exemplars $\mathcal{P}_{24}(x_0x_3x_6x_1)$ and $\mathcal{P}_{24}(x_0x_1x_2x_{19})$; see Lemma \ref{Lemma:FreeProduct} and Section \ref{Section:Conclusion}.  

\bigskip

Our main results regarding asphericity and finiteness are Theorems \ref{Theorem:AsphericityClassification} and \ref{Theorem:FinitenessClassification}. Note that the finiteness classification is complete, with no unresolved cases remaining. 
\begin{maintheorem}[Asphericity]\label{Theorem:AsphericityClassification} If $\mathcal{P} = \mathcal{P}_n(x_0x_jx_kx_l)$ is not of type $(U^\ast)$, then $\mathcal{P}$ is combinatorially aspherical if and only if at least one of the following holds:
\begin{enumerate}
\item[(a)] (A) and  (C) are both true, 
\item[(b)] (B) and (C) are both false, or
\item[(c)] (B) is true, (A) and (C) are false, and $\mathcal{P}$ is not of type $(I^\ast)$.
\end{enumerate}
If $\mathcal{P}_n(x_0x_jx_kx_l)$ is combinatorially aspherical and either $k \not \equiv 0$ or $j \not \equiv l$ mod $n$, then $G_n(x_0x_jx_kx_l)$ is a torsion-free infinite group of geometric dimension at most two.
\end{maintheorem}

\begin{maintheorem}[Finiteness]\label{Theorem:FinitenessClassification} Given $(n,j,k,l)$, the cyclically presented group $G = G_n(x_0x_jx_kx_l)$ is finite if and only if at least one of the following holds:
\begin{enumerate}
\item[(a)] $n=1$;
\item[(b)] $\gcd(n,2k) =1$ and either
\begin{itemize}
\item $l \equiv j+k$ mod $n$ and $\gcd(n,j) = 1$, or
\item $l \equiv j-k$ mod $n$ and $\gcd(n,l) = 1$;
\end{itemize}
\item[(c)] (A) and (B) are true, (C) is false, and the secondary divisor $\gamma = 1$;
\item[(d)] The presentation $\mathcal{P}_n(x_0x_jx_kx_l)$ is of type (I5), (I6$^{\,\prime}$), or (I6$^{\,\prime\prime}$).
\end{enumerate}
In parts (a) and (b), $G \cong \Z_4$. In part (c), $G$ is solvable. In type (I5), $G$ is metacyclic and non-nilpotent of order $220$. The types (I6$^{\,\prime}$) and (I6$^{\,\prime\prime}$) lead to nonisomorphic nonsolvable groups $G$ of order $4\,088\,448 = 2^7\cdot 3^3 \cdot 7 \cdot 13^2$, each containing the simple group $\mathrm{PSL}(3,3)$.
\end{maintheorem}

These classifications are supported by simultaneous consideration of \textbf{shift dynamics}, which concern the action of the cyclic group $\Z_n$ by automorphisms on (the non-identity elements of) $G = G_n(w)$. For arbitrary $n$ and $w$, shift dynamics encompass at least these two general problems.
\begin{description}
\item[Fixed Points:] For which $n$ and $w$ does the shift on $G_n(w)$ have a nonidentity fixed point?
\item[Freeness:] For which $n$ and $w$ does the shift determine a free $\Z_n$-action on the nonidentity elements of $G_n(w)$? 
\end{description}

Asphericity and freeness are related by a general result from \cite{BShift}. A cyclic presentation $\mathcal{P}_n(w)$ is \textbf{orientable} if $w$ is not a cyclic permutation of the inverse of any of its shifts. In the general case, a cyclic presentation $\mathcal{P}_n(w)$ fails to be orientiable if and only if $n = 2m$ is even and $w = u\theta_F^m(u)^{-1}$ for some word $u$ \cite[Lemma 3.6]{BShift} (in which case $u \in \mathrm{Fix}(\theta_G^m)$). Note that any cyclic presentation defined by a positive word is orientable. 

\begin{Theorem}\emph{\textbf{(\cite[Theorem A]{BShift})}}\label{Theorem:Free} For all positive $n$ and arbitrary $w$, if the cyclic presentation $\mathcal{P}_n(w)$ is orientable and combinatorially aspherical, then $\Z_n$ acts freely via the shift on the non-identity elements of the cyclically presented group $G_n(w)$. 
\end{Theorem}
As a consequence, if a combinatorially aspherical cyclic presentation $\mathcal{P}_n(w)$ is orientable and defines a nontrivial group $G = G_n(w)$, then the shift automorphism $\theta_G$ has order $n$ (compare \cite{JExt}) and in fact determines an element of order $n$ in the \textit{outer} automorphism group $\mathrm{Out}(G)$. This is because every inner automorphism of a nontrivial group has a nonidentity fixed point. We note that there are just a few known examples of orientable cyclic presentations that are not combinatorially aspherical and yet have free shift action. Examples include the Fibonacci groups $F(2,n) = G_n(x_0x_1x_2^{-1})$ for $n = 5,7$, which are finite cyclic having orders  $11$ and $29$ respectively \cite{TRevisit}, but these appear to be small exceptions rather than typical cases. 

%Cyclic presentations with positive relators are orientable and it appears that the links between Asphericity and Freeness, as well as between Finiteness and Fixed Points, are especially close.  or positive relators of length In contrast, there are the following general results for positive relators of length three. Note that positive words are orientable.
%
%\begin{Theorem}[\cite{BShift}] \label{Theorem:Length3} Let $n$ be a positive integer and let $w = x_0x_kx_l$ be a positive word of length three. 
%\begin{enumerate}
%\item[(a)] The cyclic presentation $\mathcal{P}_n(w)$ is combinatorially aspherical if and only if $\Z_n$ acts freely on the non-identity elements of the cyclically presented group $G_n(w)$. 
%\item[(b)] The cyclically presented group $G_n(w)$ is finite iff and only if the shift on $G_n(w)$ has a nonidentity fixed point.
%\end{enumerate}
%\end{Theorem}
%
In \cite[Theorems B,C]{BShift} it was shown that if $w = x_0x_kx_l$ is a positive word of length three, then the relationships between finiteness, asphericity, and shift dynamics are sharp:
\begin{itemize}\em
\item $\mathcal{P}_n(x_0x_kx_l)$ is combinatorially aspherical if and only if $\Z_n$ acts freely via the shift on the nonidentity elements of $G_n(x_0x_kx_l)$;
\item $G = G_n(x_0x_kx_l)$ is finite if and only if the shift $\theta_G$ has a nonidentity fixed point.
\end{itemize}
One of the principal aims of this article is to extend this connection to the case of length four relators. The relation between asphericity and freeness is not as tidy due to the unresolved cases in the asphericity classification.

\begin{maintheorem}\label{Theorem:AsphericalFree} If $\mathcal{P} = \mathcal{P}_n(x_0x_jx_kx_l)$ is not of type $(U^\ast)$, then $\mathcal{P}$ is combinatorially aspherical if and only if the cyclic group $\Z_n$ acts freely on the non-identity elements of $G_n(x_0x_jx_kx_l)$. 
\end{maintheorem}
All available evidence suggests that the identification between asphericity and free action will persist once the asphericity status is resolved for type (U$^\ast$).

When combined with Theorem \ref{Theorem:FinitenessClassification}, the next result explicitly determines those cyclic presentations $\mathcal{P}_n(x_0x_jx_kx_l)$ for which the shift has a nonidentity fixed point.
\begin{maintheorem}\label{Theorem:FiniteFix} The shift on $G = G_n(x_0x_jx_kx_l)$ has a non-identity fixed point if and only if $G$ is finite or else (C) is true and $\gcd(n,2k) = 1$.
\end{maintheorem}
Thus it can happen that $G = G_n(x_0x_jx_kx_l)$ is infinite and $\mathrm{Fix}(\theta_G) \neq 1$. By analyzing the infinite cases in detail, we are able to show that nonidentity fixed points for the shift always derive from finite cyclically presented subgroups.

\begin{maintheorem}\label{Theorem:Finite=Fix} The shift on $G = G_n(x_0x_jx_kx_l)$ has a nonidentity fixed point if and only if $G$ possesses a finite cyclically presented subgroup of the form $H = G_n(v)$ where $\theta_G|_H = \theta_H$ and $\mathrm{Fix}(\theta_G) = \mathrm{Fix}(\theta_H) \neq 1$. %In particular, $\mathrm{Fix}(\theta_G)$ is finite for all $(n,j,k,l)$. 
\end{maintheorem}
%In Section \ref{Section:Conclusion} we show that if $G = G_n(x_0x_jx_kx_l)$ is infinite and $\mathrm{Fix}(\theta_G) \neq 1$, then $w = x_0x_jx_kx_l$ (or its cyclic permutation $x_jx_kx_lx_0$) is expressible as a \textbf{composite} word $w = v \circ u$ where $u$ and $v$ are length two positive words. 

Combining with the results of \cite{BShift}, we conclude the following.

\begin{Corollary} If $G = G_n(w)$ where $n$ is a positive integer and $w$ is a positive word of length at most four, then the fixed point subgroup $\mathrm{Fix}(\theta_G)$ is finite. Moreover, if $G$ is finite then $\mathrm{Fix}(\theta_G)$ is nontrivial.
\end{Corollary}

We have already outlined the contents of the preliminary sections \ref{Section:Asphericity}-\ref{Section:ABC}, so we conclude this introduction with a comment on the structure of the remainder of the article. For each combination of the conditions (A), (B), (C), we treat finiteness, asphericity, and shift dynamics simultaneously. In this way we take advantage of the intrinsic connections between these concepts and avoid needless repetition of common arguments. Section \ref{Section:SmallCancellation} treats four of the eight possible combinations of (A), (B), (C), all of which lead to small cancellation conditions. Section \ref{Section:FxT} deals with the case where (A) is false and (C) is true (FxT). Section \ref{Section:TTF} deals with the case where (A) and (B) are true and (C) is false (TTF). Section \ref{Section:FTF} deals with the most interesting case, where (B) is true and both (A) and (C) are false (FTF). A concluding section wraps up the details of Theorems \ref{Theorem:AsphericityClassification}-\ref{Theorem:Finite=Fix}.

The results of this article are drawn from the PhD dissertation of the second named author. We are grateful to Gerald Williams for helpful comments and for granting permission to discuss his computations regarding the group $L$ with presentation $\pres{t,u}{t^6, u^3t^3ut}$, which plays key roles in Lemma \ref{Lemma:E6} and in Theorem \ref{Theorem:FinitenessClassification}.

\section{Combinatorial and Relative Asphericity}\label{Section:Asphericity}

%Our immediate goal is to relate combinatorial asphericity of a cyclic presentation $\mathcal{P}_n(w)$ to the asphericity of a certain relative presentation for the shift extension $E_n(w)  = G_n(w) \rtimes_\theta \Z_n$. The organization of this section follows that of \cite[Sections 2,3]{BShift}. We begin by discussing the concept of combinatorial asphericity.

Given a presentation $\mathcal{P} = \pres{\mathbf{x}}{\mathbf{r}}$ for a group $G$, the free group $F = F(\mathbf{x})$ with basis $\mathbf{x}$ acts by automorphisms on the free group $\mathbb{F} = F(F\times \mathbf{r})$ with basis $F \times \mathbf{r}$ by $v \cdot (u,r) = (vu,r)$. The homomorphism $\partial: \mathbb{F} \ra F$ given by $\partial(u,r) = uru^{-1}$ is $F$-equivariant where $F$ acts on itself by (left) conjugation. The image of $\partial$ is the normal closure of $\mathbf{r}$ in $F$ so the cokernel of $\partial$ is isomorphic to $G$. Elements of the kernel $\mathbb{I} = \ker \partial$ are called \textbf{identities} (or \textbf{identity sequences}) for $\mathcal{P}$. Then we have the exact sequence 
\begin{equation}\label{Equation:FundSeq1}
1 \ra \mathbb{I} \ra \mathbb{F} \stackrel{\partial}{\ra} F \ra G \ra 1
\end{equation}
for the presentation $\mathcal{P}$. Among the identities are the so-called \textbf{Peiffer identities} having the form $(u,r)^\epsilon(v,s)^\delta(u,r)^{-\epsilon}(ur^\epsilon u^{-1}v,s)^{-\delta}$ where $u,v \in F$, $r,s \in {\bf r}$, and $\delta, \epsilon = \pm 1$. If $\mathbb{P}$ denotes the normal closure of the Peiffer identities in $\mathbb{F}$, then the $F$-action on $\mathbb{F}$ descends to a $G$-action on the abelian quotient group $\mathbb{I}/\mathbb{P}$ and there is the exact \textbf{fundamental sequence}
\begin{equation}\label{Equation:FundSeq}
0 \ra \mathbb{I}/\mathbb{P} \ra \mathbb{F}/\mathbb{P} \stackrel{\partial}{\ra} F \ra G \ra 1
\end{equation}
that provides a combinatorial description for the long exact homotopy sequence of the skeleton pair $(K,K^1)$ where $K$ is the two-complex modeled on the presentation $\mathcal{P}$. In particular $\pi_2 K \cong \mathbb{I}/\mathbb{P}$ as $\Z G$-modules. See \cite[Section 1]{CCH}, \cite[Section 2]{AJSFrame}, or \cite{AJSAlgTop} for details. 

As in \cite[Proposition 1.4]{CCH}, a presentation $\mathcal{P} = \pres{\bf x}{\bf r}$ for a group $G$ is \textbf{combinatorially aspherical} if the homotopy module $\pi_2 K \cong \mathbb{I}/\mathbb{P}$ is generated as a $\Z G$-module by the classes of identities that have length two in the free group $\mathbb{F}$ with basis $F(\mathbf{x}) \times \mathbf{r}$. The two-dimensional cellular model $K$ of the presentation $\mathcal{P}$ is aspherical, in the sense that $\pi_2K = 0$, if and only if the presentation $\mathcal{P}$ is combinatorially aspherical and no relator of  $\mathcal{P}$ is a proper power, nor is conjugate to any another relator or its inverse. See \cite{CCH}, \cite[Section 3]{BShift}, or the survey \cite{PIdentities}, which contains a thorough discussion of identity sequences and their relation to \textbf{pictures} under the blanket assumption that no relator is conjugate to any another relator or its inverse. From the point of view of spherical pictures (or dually, spherical diagrams), an arbitrarily given presentation is combinatorially aspherical if and only if its second homotopy module is generated by the classes of spherical pictures that contain exactly two relator discs. See \cite[Section 2.3]{PIdentities}.

Presentations satisfying any of the C(p)-T(q) small cancellation conditions are combinatorially aspherical if $1/p+1/q\leq1/2$  \cite{CCH,Lyndon66}. In fact they satisfy even stronger conditions such as diagrammatic asphericity \cite{CH}. Here we shall be interested in the C(4)-T(4) condition. Recall that a word is a \textbf{piece} for a presentation if it occurs as a common initial subword of two distinct cyclic permutations of the relators or their inverses.

\begin{Lemma}\label{Lemma:SmCancCrit} The presentation $\mathcal{P}_n(x_0x_jx_kx_l)$ satisfies the C(4)-T(4) small cancellation condition if and only if no length two cyclic subword of $w = x_0x_jx_kx_l$ is a piece.
\end{Lemma}

\proof Since the relators of $\mathcal{P}_n(x_0x_jx_kx_l)$ have length four, the C(4) condition is satisfied if and only if each piece has length one. On the other hand, if all pieces have length one, then the T(4) condition is also satisfied because the relators are all positive words, see \cite{HPV}. The cyclic symmetry in the presentation implies that it suffices to verify that no length two cyclic subword of $w=x_0x_jx_kx_l$ is a piece. \qed

\bigskip

A criterion for combinatorial asphericity that is specific to cyclic presentations involves the concept of an aspherical relative presentation. Given a cyclic presentation $\mathcal{P}_n(w)$, the shift action by the cyclic group $\Z_n$ on $G = G_n(w)$ determines the \textbf{shift extension}

$$
E = E_n(w) = G \rtimes_{\theta_G} \Z_n,
$$
which admits a two-generator two-relator presentation of the form $\pres{a,x}{a^n, W}$ where $W$ is obtained from $w$ via the substitutions $x_i = a^ixa^{-i}$. The shift $\theta_{G}$ then arises via conjugation by $a$ in the shift extension and $G$ is recovered as the kernel of a retraction $\nu:E \ra \Z_n$ that satisfies $\nu(a) = a$ and $\nu(x) = 1$. As in \cite[Theorem 2.3]{BShift} if $E$ is any group with a presentation of the form $E = \pres{a,x}{a^n,W}$ and $\nu: E \ra \pres{a}{a^n}$ is a retraction onto the cyclic group of order $n$ generated by $a$, then the kernel is cyclically presented: $\ker \nu \cong G_n(\rho^\nu(W))$ where $\rho^\nu$ is a Reidemeister-Schreier rewriting process that depends on the retraction $\nu$. For example if $\nu(x) = a^0 = 1$, then $\rho^\nu(xa^jxa^{k-j}xa^{l-k}xa^{-l}) = x_0x_jx_kx_l$. Given $E \cong \pres{a,x}{a^n,W}$, there can be multiple retractions and the cyclically presented kernels arising from different retractions need not be isomorphic. However these retraction kernels are obviously commensurable in $E$ and \cite[Lemma 2.2]{BShift} implies that they have identical shift dynamics.

%\begin{equation}\label{Equation:CyclicPres}
%\mathcal{P}_n(x_0x_jx_kx_l) = \pres{x_0, \ldots, x_{n-1}}{x_ix_{i+j}x_{i+k}x_{i+l}, 0 \leq i < n}
%\end{equation}
%for $G = G_n(x_0x_jx_kx_l)$. The group $G$ is nontrivial, having the cyclic group $\pres{x}{x^4}$ of order four as a homomorphic image by mapping each $x_i \mapsto x$. The shift extension is $E = E_n(x_0x_jx_kx_l) = G \rtimes_{\theta_G} \Z_n$ with ordinary and relative presentations
%\begin{equation}\label{Equation:RelativePres}
% E \cong \pres{a,x}{a^n, xa^jxa^{k-j}xa^{l-k}xa^{-l}} \cong \pres{\Z_n, x}{xa^jxa^{k-j}xa^{l-k}xa^{-l}}
%\end{equation}
%where the coefficient group $\Z_n \cong \pres{a}{a^n}$ is cyclic of order $n$, generated by $a$. Viewing $G$ as the kernel of a retraction $\nu:E \ra \Z_n$ that trivializes $x$, the shift $\theta_G \in \mathrm{Aut}(G)$ on the normal subgroup $G \unlhd E$ is given by conjugation by $a$: $\theta_G(g) = aga^{-1}$ for all $g \in G$. The word $w = x_0x_jx_kx_l$ is recovered via the rewriting $w = \rho^\nu(W)$ where $W = xa^jxa^{k-j}xa^{l-k}xa^{-l}$. 
%
%This means that the fixed point set of for any power of the shift $\theta_G^d$ is expressible in terms of a centralizer in $E = E_n(w)$:
%
%$$
%\mathrm{Fix}(\theta_G^d) = G \cap \mathrm{Cent}_E(a^d).
%$$
%The cyclically presented group $G_n(w)$ is recovered from $E_n(w)$ as the kernel of the retraction $E_n(w) \ra \Z_n$ onto the cyclic subgroup generated by $a$ that sends $x$ to the identity. %
%

Given any presentation $\pres{a,x}{a^n,W}$ for a group $E$, we can view the word $W$ in the free product $\Z_n \ast \sgp{x} \cong \pres{a,x}{a^n}$ to obtain a one-relator \textbf{relative presentation} $\pres{\Z_n,x}{W}$ with \textbf{coefficient group} $\Z_n \cong \pres{a}{a^n}$. A \textbf{cellular model} $M$ and a relativized concept of asphericity were described in \cite[Section 4]{BP}. Starting with an Eilenberg-Maclane complex of type $K(\Z_n,1)$, we first form the complex $K(\Z_n,1) \vee S^1_x \cup c^2_W$, which is the model $M$ if $W$ is not a proper power when viewed in the free product $\Z_n \ast \sgp{x}$. Otherwise, if $W = \dot{W}^e$ in the free product where $e > 1$, then we add cells in dimensions three and higher to obtain the model $M$. The homotopy type of the model $M$ is determined by the relative presentation $\pres{\Z_n,x}{W}$. See \cite[Section 3]{BShift} or \cite[Section 4]{BP} for details. The relative presentation $\pres{\Z_n, x}{W}$ is \textbf{aspherical} if the relative homotopy group $\pi_2(M,K(\Z_n,1))$ is trivial. Equivalently, the inclusion $K(\Z_n,1) \ra M$ induces an injective homomorphism on fundamental group and the absolute homotopy group $\pi_2 M$ is trivial. The point is that if $\pres{\Z_n,x}{W}$ is an aspherical relative presentation for a group $E$, then $\Z_n$ embeds in $E$ and $M$ is a $K(E,1)$ complex. In the situation where $E = \pres{a,x}{a^n,W}$ is the shift extension for a cyclically presented group $G_n(w)$, then $\Z_n$ is a retract of $E = \pi_1M$ and so the embedding condition on fundamental groups is automatic; thus the relative presentation $\pres{\Z_n,x}{W}$ is aspherical if and only $\pi_2M = 0$. 

If a group $E$ with relative presentation $\pres{\Z_n,x}{W}$ admits a retraction $\nu: E \ra \Z_n$ then there is the following connection between asphericity of the relative presentation $\pres{\Z_n,x}{W}$ and combinatorial asphericity of the corresponding cyclic presentation $\mathcal{P}_n(\rho^\nu(W))$. 

\begin{Theorem}\emph{\textbf{(\cite[Theorem 4.1]{BShift})}}\label{Theorem:AsphTransf} Let $M$ be the cellular model of a relative presentation $\pres{\Z_n,x}{W}$ for a group $E$. Suppose that $\nu:E \ra \Z_n$ is a retraction and let $w = \rho^\nu(W)$.
\begin{enumerate}
  \item[(a)] If $\pi_2 M = 0$, then $\mathcal{P}_n(w)$ is combinatorially aspherical.
  \item[(b)] If $\mathcal{P}_n(w)$ is orientable and combinatorially aspherical, then $\pi_2M = 0$.
\end{enumerate}
\end{Theorem}

%We will make use of the following algebraic consequence of relative asphericity, here specialized to the current context from \cite{BP}.
%
%\begin{Theorem}\emph{\textbf{(\cite[(0.4)]{BP})}}\label{Theorem:FiniteSubgroup} Suppose that $\pres{\Z_n,x}{W}$ is an aspherical relative presentation for a group $E$ and that $W = \dot{W}^e$ in the free product $\Z_n \ast \sgp{x}$ where the exponent $e$ is maximal. Then the element $\dot{W}$ determines an element of order $e$ in $E$ and every finite subgroup of $E$ is conjugate to a subgroup of either $\Z_n$ or $\sgp{\dot{W}}$. 
%\end{Theorem}

\section{Isomorphisms}\label{Section:DU}

Classifications involving cyclic presentations are complicated by the fact that there are many isomorphisms that must be accounted for, not all of which are obvious or easily catalogued. See \cite[Section 2]{CRS05} and \cite[Lemma 2.1]{EdjvetWilliams}, for example.

For a fixed position integer $n$, let $\Phi_n$ denote the set of all positive words of length four in the generators $x_0, \ldots, x_{n-1}$ for the free group $F$ of rank $n$. Thus every element $w \in \Phi_n$ has the form $w = x_ix_jx_kx_l$ where integers $i,j,k,l$ are defined modulo $n$. We consider the following bijective transformations on $\Phi_n$:
\begin{eqnarray*}
\sigma(x_ix_jx_kx_l) &=& x_ix_lx_kx_j\\
\tau(x_ix_jx_kx_l) &=& x_lx_ix_jx_k\\
\theta_F(x_ix_jx_kx_l) &=& x_{i+1}x_{j+1}x_{k+1}x_{l+1}\\
u(x_ix_jx_kx_l) &=& x_{ui}x_{uj}x_{uk}x_{ul}.
\end{eqnarray*}
Here the element $u$ is taken from the multiplicative group $\Z_n^\ast$ of units modulo $n$. All subscripts are interpreted modulo $n$. 

\begin{Lemma}\label{Lemma:DU} Within the group $\mathrm{Sym}(\Phi_n)$ of permutations of $\Phi_n$, the transformations $\sigma, \tau, \theta_F$ and $u \in \Z_n^\ast$ generate a subgroup that is a homomorphic image of $\Gamma_n = D_4 \times (\Z_n \rtimes \Z_n^\ast)$, where $D_4$ is the dihedral group of order $8$ generated by $\sigma$ and $\tau$ and the semi-direct product involves the natural action of $\Z_n^\ast \cong \mathrm{Aut}(\Z_n)$ acting on the additive group $\Z_n$ by multiplication.
\end{Lemma}

\proof Working in $\mathrm{Sym}(\Phi_n)$, one checks that $\sigma^2 = \tau^4 = (\sigma \tau)^2 = \theta_F^n = 1$ and $u\theta_F = \theta_F^u u$ for all $u \in \Z_n^\ast$. Moreoever the transformations $\sigma, \tau$ commute pairwise with $\theta_F, u
$. \qed

\bigskip

\noindent The orbit of $w = x_ix_jx_kx_l \in \Phi_n$ under the action of $\Gamma_n$ consists of all words that arise as cyclic permutations of shifts of words of the form $u(w) = x_{ui}x_{uj}x_{uk}x_{ul}$ or $\sigma(u(w)) = x_{ui}x_{ul}x_{uk}x_{uj}$ where $u \in \Z_n^\ast$. 

\begin{Theorem}\label{Theorem:DUI} For each $c \in \Gamma_n = D_4 \times (\Z_n \rtimes \Z_n^\ast)$ and $w = x_ix_jx_kx_l \in \Phi_n$, there is a group isomorphism $c_w: G_n(w) \ra G_n(c(w))$ such that 
\begin{equation}\label{Equation:ShiftEq}
c_w \circ \theta_{G_n(w)} = \theta_{G_n(c(w))}^{\psi(c)} \circ c_w
\end{equation}
where $\psi: \Gamma_n \ra \Z_n^\ast$ is the natural projection. Moreover $\mathcal{P}_n(w)$ is combinatorially aspherical if and only if $\mathcal{P}_n(c(w))$ is combinatorially aspherical.
\end{Theorem}

\proof Let $F = F(\mathbf{x})$ be the free group with basis $\mathbf{x} = \{x_p\ : 0 \leq p < n\}$. Given $c \in \Gamma_n$ and $w = x_ix_jx_kx_l  \in \Phi_n$, we denote the relator sets for the presentations $\mathcal{P}_n(w)$ and $\mathcal{P}_n(c(w))$ by $$\mathbf{r}_w = \{\theta_F^q(w)\ :\ 0 \leq q < n\} = \{x_{i+q}x_{j+q}x_{k+q}x_{l+q}\ :\ 0 \leq q < n\}$$ and $$c(\mathbf{r}_w) = \{\theta_F^q(c(w))\ :\ 0 \leq q < n\}$$ respectively. We first construct homomorphisms $c_w: F \ra F$ and $\hat{c}_w: \mathbb{F}_w \ra \mathbb{F}_{c(w)}$ where $\mathbb{F}_w = F(F \times \mathbf{r}_w)$ and similarly $\mathbb{F}_{c(w)} = F(F \times c(\mathbf{r}_w))$. For $c = \sigma, \tau, \theta_F$ or $u \in \Z_n^\ast$ and $w \in \Phi_n$, the assignments
\begin{eqnarray*}
\sigma_w(x_p) =  x_p^{-1}; && \hat{\sigma}_w(v,\theta_F^q(w)) = (\sigma_w(v)\theta_F^q(x_i)^{-1},\theta_F^q\sigma(w))^{-1}\\
\tau_w = 1_F; && \hat{\tau}_w(v,\theta_F^q(w)) = (v\theta_F^q(x_l)^{-1},\theta_F^q\tau(w))\\
(\theta_F)_w = \theta_F;&& \widehat{(\theta_F)}_w(v,\theta_F^q(w)) = (\theta_F(v),\theta_F^{q+1}(w))\\
u_w(x_p) = x_{up};&& \hat{u}_w(v,\theta_F^q(w)) = (u_w(v), \theta_F^{uq}(u(w))).
\end{eqnarray*}
define length-preserving homomorphisms $c_w: F \ra F$ and $\hat{c}_w: \mathbb{F}_w \ra \mathbb{F}_{c(w)}$. These homomorphisms are constructed to satisfy $\partial \circ \hat{c}_w = c_w \circ \partial$ for $c = \sigma, \tau, \theta_F, u \in \Gamma_n = D_4 \times (\Z_n \rtimes \Z_n^\ast)$. Viewed in $\mathrm{Aut}(F)$, the isomorphisms $\sigma_w, \tau_w, \theta_F, u_w$ are compatible with the relations of $\Gamma_n = D_4 \times (\Z_n \rtimes \Z_n^\ast)$. Specifically, one notes that
$$
\sigma_w^2 = \tau_w^4 = \theta_F^n = 1, t_w \circ u_w = (tu)_w,\ \ (t,u \in \Z_n^\ast)
$$
and that the homomorphisms $\sigma_w, \tau_w$ commute pairwise with $\theta_F$ and $u_w\ \ (u \in \Z_n^\ast)$. Since $\Gamma_n$ acts on $\Phi_n$, this shows that $\Gamma_n$ acts on $F$ by automorphisms and in fact this action is independent of $w$. (We retain the notation $c_w$ in order to distinguish between $c \in \Gamma_n$ and $c_w \in \mathrm{Aut}(F)$.) For each $c \in \Gamma_n$ we now have an automorphism $c_w \in \mathrm{Aut}(F)$ that is length-preserving and which is compatible with the boundary maps $\partial$ from the fundamental sequences (\ref{Equation:FundSeq1}) for the presentations $\mathcal{P}_n(w)$ and $\mathcal{P}_n(c(w))$.%\footnote{This can be viewed as constructing specific homotopy equivalences on the skeleton pairs $(K,K^1)$ for the corresponding cellular models. Yet another perspective would be to view these equivalences at the level of the models of the corresponding relative presentations.}
$$
\begin{CD} 
{\mathbb{F}_w = F(F \times \mathbf{r}_w)} @> {\partial} >> {F} \\ 
@ V{\hat{c}_w}VV @ VV{c_w} V  \\ 
{\mathbb{F}_{c(w)} = F(F \times c(\mathbf{r}_w))} @> {\partial} >> {F}\\ 
\end{CD} 
$$
It follows that for any $c \in \Gamma_n$, the automorphism $c_w \in \mathrm{Aut}(F)$ induces an isomorphism $G_n(w) \ra G_n(c(w))$, which we also denote by $c_w$. As for shift equivariance, we have

\begin{eqnarray*}
\sigma_w \circ \theta_{G_n(w)} &=& \theta_{G_n(\sigma(w))} \circ \sigma_w\\
\tau_w \circ \theta_{G_n(w)} &=& \theta_{G_n(\tau(w))} \circ \tau_w\\
u_w \circ \theta_{G_n(w)} &=& \theta^u_{G_n(u(w))} \circ u_w.
\end{eqnarray*}
Moreover $G_n(\theta_F(w)) = G_n(w)$ and $(\theta_F)_w = \theta_{G_n(w)}$. This verifies equivariance as in (\ref{Equation:ShiftEq}). 

It remains to verify that combinatorial asphericity is preserved under the action of $\Gamma_n$ on the set of presentations of the form $\mathcal{P}_n(w)$, $w \in \Phi_n$. Note that the homomorphisms $\hat{\sigma}_w, \hat{\tau}_w$ depend on $w$. It turns out that the homomorphisms $\hat{\sigma}_w, \hat{\tau}_w, \widehat{(\theta_F)}_w, \hat{u}_w$ are compatible with some of the relations of $\Gamma_n$ but not all. To illustrate that $\hat{\sigma}_w^2 = 1_{\mathbb{F}_w}$, we use the fact that $\sigma^2 = 1_{\Phi_n}$, that $\sigma_w$ commutes with $\theta_F$ in $\mathrm{Aut}(F)$, and that the first letter of $w=x_ix_jx_kx_l$ is the same as for $\sigma(w) = x_ix_lx_kx_j$, we have

\begin{eqnarray*}
\hat{\sigma}_w^2(v,\theta_F^q(w)) &=& \hat{\sigma}_w(\sigma_w(v)\theta_F^q(x_i)^{-1},\theta_F^q\sigma(w))^{-1}\\
&=&(\sigma_w\left(\sigma_w(v)\theta_F^q(x_i)^{-1}\right) \cdot \theta_F^q(x_i)^{-1}, \theta_F^q\sigma^2(w))\\
&=&(v \theta_F^q(\sigma_w(x_i)^{-1}) \cdot \theta_F^q(x_i)^{-1}, \theta_F^q(w))\\
&=& (v \theta_F^q(x_i) \cdot \theta_F^q(x_i)^{-1}, \theta_F^q(w)) = (v,  \theta_F^q(w)).
\end{eqnarray*}
However $\hat{\tau}_w^4: \mathbb{F}_w \ra \mathbb{F}_w$ is not the identity. For example, 

\begin{eqnarray*}
\hat{\tau}_w^4(1,x_ix_jx_kx_l) &=& \hat{\tau}_w^3(x_l^{-1},x_lx_ix_jx_k)\\
&=&\hat{\tau}_w^2(x_l^{-1}x_k^{-1},x_kx_lx_ix_j)\\
&=&\hat{\tau}_w(x_l^{-1}x_k^{-1}x_j^{-1},x_jx_kx_lx_i)\\
&=&(x_l^{-1}x_k^{-1}x_j^{-1}x_i^{-1},x_ix_jx_kx_l)
\end{eqnarray*}
Nevertheless, the fact that $x_l^{-1}k_k^{-1}kx_j^{-1}x_i^{-1} = 1$ in the group $G_n(w) = G_n(x_ix_jx_kx_l)$ implies that $(1,x_ix_jx_kx_l) \equiv (x_l^{-1}x_k^{-1}x_j^{-1}x_i^{-1},x_ix_jx_kx_l)$ modulo the group $\mathbb{P}_w$ of Peiffer identities for the presentation $\mathcal{P}_n(w)$. See \cite[Lemma 2.4]{AJSAlgTop}. In fact one can verify that the homomorphisms $\hat{\sigma}_w, \hat{\tau}_w, \widehat{(\theta_F)}_w, \hat{u}_w$ are compatible with all of the relations of $\Gamma_n$ if we work modulo Peiffer identities. For each $c \in \Gamma_n$ and $w \in \Phi_n$, we obtain a commutative diagram that compares the fundamental sequences (\ref{Equation:FundSeq}) for the presentations $\mathcal{P}_n(w)$ and $\mathcal{P}_n(c(w))$. The vertical arrows all represent isomorphisms. 
%
%Thus our claim is that if we work modulo Peiffer identities, then for all $w \in \Phi_n$ and $c = \sigma, \tau, \theta_F, u \in \Gamma_n$, we do obtain homomorphisms $\bar{c}_w: F(F \times \mathbf{r})/\mathbb{P}_w \ra F(F \times c(\mathbf{r}))/\mathbb{P}_{c(w)}$ that are compatible with the relations of $\Gamma_n$. Here we let $\mathbb{P}_w$ stand for the group of Peiffer identities for the presentation $\mathcal{P}_n()w)$ and similarly for $\mathbb{P}_{c(w)}$. For example we have $\bar{\tau}_w^4 = 1_{\mathbb{F}/\mathbb{P}}$. From this it follows that for all $c \in \Gamma_n$, we obtain homomorphisms $c_w: F \ra F$ and $\bar{c}_w: F(F \times \mathbf{r})/\mathbb{P}_w \ra F(F \times c(\mathbf{r}))/\mathbb{P}_{c(w)}$ that make the following diagram commute.
$$
\begin{CD} 
{0}@>>> {\mathbb{I}_w/\mathbb{P}_w} @>>>{\mathbb{F}_w/\mathbb{P}_w} @> {\partial} >> {F} @>>>G_n(w) @>>> 1 \\ 
@.@ VVV @V{\hat{c}_w}VV @ VV{c_w} V @VVV \\ 
{0}@>>> {\mathbb{I}_{c(w)}/\mathbb{P}_{c(w)}} @>>>{\mathbb{F}_{c(w)}/\mathbb{P}_{c(w)}} @> {\partial} >> {F} @>>>G_n(c(w)) @>>> 1 \\ 
\end{CD} 
$$
%
%
%\begin{eqnarray*}
%\sigma_w^2 = 1_F; && \hat{\sigma}_w^2 = 1_{F(F \times \mathbf{r}_w)}\\
%\tau_w^4 = 1_F; && \hat{\tau}_w(v,\theta^q(w)) = (\tau_w(v)\theta_F^q(x_l)^{-1},\theta_F^q\tau(w))\\
%(\theta_F)_w^n = 1_F;&& \widehat{(\theta_F)}_w(v,\theta^q(w)) = (\theta_F(v),\theta_F^{q+1}(w))\\
%t_w \circ u_w = (tu)_w;&& \hat{u}_w(v,\theta^q(w)) = (u_w(v), \theta_F^{uq}(u(w))).
%\end{eqnarray*}
%
%
%
%preserve the appropriate relations and so define homomorphisms $c_w: G_n(w) \ra G_n(c(w))$. These homomorphisms are compatible with the defining relations for $\Gamma_n$, as for example $\sigma_w^2 = 1_{G_n(w)}$ where $\sigma^2(w) = w$ and $u_w \circ (\theta_F)_w = (\theta_F)_w^u \circ u_w$ as homomomrphisms $G_n(w) \ra G_n(\theta_F(u(w)) = G_n(\theta_F^u(u(w)))$. Thus for arbitrary $c \in \Gamma_n$, isomorphisms $c_w$ are defined via composition because $\Gamma_n$ acts on $\Phi_n$ as in Lemma \ref{Lemma:DU}. As for shift equivariance as in (\ref{Equation:ShiftEq}), we have
%
%\begin{eqnarray*}
%\sigma_w \circ \theta_{G_n(w)} &=& \theta_{G_n(\sigma(w))} \circ \sigma_w\\
%\tau_w \circ \theta_{G_n(w)} &=& \theta_{G_n(\tau(w))} \circ \tau_w\\
%u_w \circ \theta_{G_n(w)} &=& \theta^u_{G_n(u(w))} \circ u_w.
%\end{eqnarray*}
%Moreover $G_n(\theta_F(w)) = G_n(w)$ and $(\theta_F)_w = \theta_{G_n(w)}$. 
It remains to note that the isomorphism $\mathbb{I}_w/\mathbb{P}_w \ra \mathbb{I}_{c(w)}/\mathbb{P}_{c(w)}$ is length preserving and so it follows that $\mathbb{I}_w/\mathbb{P}_w$ is generated by classes of length two identity sequences if and only if the same is true for $\mathbb{I}_{c(w)}/\mathbb{P}_{c(w)}$.  \qed

\paragraph{Remark} Given $w = x_ix_jx_kx_l \in \Phi_n$, let $M_w$ be the cellular model of the relative presentation $\pres{\Z_n,x}{a^ixa^{j-i}xa^{k-j}xa^{l-k}xa^{-l}}$ for the shift extension $E_n(w) = G_n(w) \rtimes_{\theta_{G_n(w)}} \Z_n$. Thus $M_w$ is obtained from $K(\Z_n,1) \vee S^1_x \cup c^2_w$ by attaching cells in dimensions three and up in case $w$ is a proper power in the free product $\Z_n \ast \sgp{x}$. The arguments given in the proof of Theorem \ref{Theorem:DUI} can be used to show that if we are given $c \in \Gamma_n$, then there is a homotopy equivalence $M_w \ra M_{c(w)}$ that restricts to a homotopy self-equivalence of the common subcomplex $K(\Z_n,1)$. Moreover the restricted map induces the automorphism $\psi(c) \in \mathrm{Aut}(\Z_n)$ on $\Z_n = \pi_1K(\Z_n,1)$. This leads to an isomorphism $G_n(w) \ra G_n(c(w))$ that is shift equivariant as in (\ref{Equation:ShiftEq}). That $\mathcal{P}_n(w)$ is combinatorially aspherical if and only if $\mathcal{P}_n(c(w))$ is combinatorially aspherical then follows from Theorem \ref{Theorem:AsphTransf}.

%The two-skeleton $M_w^2$ is modeled on the presentation $\pres{a,x}{a^n, a^ixa^{j-i}xa^{k-j}xa^{l-k}xa^{-l}}$ for $E_n(w)$. 

\begin{Corollary}\label{Corollary:DUDyn} If $w = x_ix_jx_kx_l \in \Phi_n$ and $c \in \Gamma_n = D_4 \times (\Z_n \rtimes \Z_n^\ast)$, then $G_n(w) \cong G_n(c(w))$ via an isomorphism that preserves shift dynamics. Thus:
\begin{itemize}
\item The shift on $G_n(w)$ has a nonidentity fixed point if and only the shift on $G_n(c(w))$ has a nonidentity fixed point; and 
\item The shift action by $\Z_n$ on the nonidentity elements of $G_n(w)$ is free if and only if the shift action by $\Z_n$ is free on the nonidentity elements of $G_n(c(w))$. 
\end{itemize}
\end{Corollary}

\proof The group isomorphisms of Theorem \ref{Theorem:DUI} are shift equivariant and so define isomorphisms of left $\Z_n$-sets, where the shift action of $\Z_n$ is twisted as necessary via an automorphism $u \in \Z_n^\ast = \mathrm{Aut}(\Z_n)$. \qed

\begin{Example} \em The isomorphisms $(\theta_F)_w$ enable the obvious fact that  in studying the presentations $\mathcal{P}_n(x_ix_jx_kx_l)$, there is no loss of generality in assuming that $i=0$. Only slightly less obvious are the isomorphisms
\begin{eqnarray*}
G_n(x_0x_jx_kx_l) &\cong& G_n(x_0x_lx_kx_j)\ \ \mbox{(via $\sigma$) and}\\
G_n(x_0x_jx_kx_l) &\cong& G_n(x_0x_{uj}x_{uk}x_{ul}) \ \ \mbox{(via $u \in \Z_n^\ast$)}.
\end{eqnarray*}
Combining the isomorphisms arising from $\sigma$ and $\tau$ we have

$$
G_n(x_0x_jx_kx_l) \cong G_n(x_0x_lx_kx_j) \cong G_n(x_lx_kx_jx_0).
$$
That is, we can read the relator backwards without changing group structure, asphericity status, or shift dynamics. 

That all of these isomorphisms preserve combinatorial asphericity can also be attributed to the fact that they arise via homotopy equivalences of the two-dimensional cellular models that are ``area-preserving" in the sense that the attaching map of each two-cell is mapped to a loop that is freely homotopic to an attaching map of a two-cell in the image. Equivalently and at the level of the free group $F$, each relator is mapped to a conjugate of a relator or its inverse in the target. In other words, the homomorphism $\mathbb{F}_w \ra \mathbb{F}_{c(w)}$ is length-preserving. These properties are witnessed in the definitions of the homomorphsims $\hat{\sigma}_w, \hat{\tau}_w, \widehat{(\theta_F)}_w, \hat{u}_w$ within the proof of Theorem \ref{Theorem:DUI}.
\end{Example}

In closing this section, we note that a similar discussion can be carried out for positive words of any given length, where positive words of length $L$ in the free group of rank $n$ are acted upon by the group $D_L \times (\Z_n \rtimes \Z_n^\ast)$. The actions give rise to group isomorphisms that preserve asphericity status and shift dynamics as in Corollary \ref{Corollary:DUDyn}. In particular, the isomorphisms of \cite[Lemma 2.1]{EdjvetWilliams} for the case $L = 3$ all arise in this way and so they too preserve shift dynamics.

\section{Divisor Criteria}\label{Section:Divisor}

To analyze cyclic presentations of the form $\mathcal{P}_n(x_0x_jx_kx_l)$, it is well known \cite{EdjvetIrred} that one can generally restrict to the case where $\gcd(n,j,k,l) = 1$. This includes considerations of (combinatorial) asphericity \cite[Theorem 4.2]{CCH}, finiteness \cite[Lemma 2.4]{CRS05}, and shift dynamics \cite[Lemma 5.1]{BShift}. The reason for this is that if $c = \gcd(n,j,k,l)$, then the cyclic presentation $\mathcal{P}_n(x_0x_jx_kx_l)$ is a disjoint union of cyclic subpresentations and $G = G_n(x_0x_jx_kx_l)$ decomposes as a free product

$$
G \cong \ast_{i=1}^c G_{n/c}(x_0x_{j/c}x_{k/c}x_{l/c})
$$
where the shift transitively permutes the free factors. In particular, if $c \neq 1$, then $G$ is infinite and the shift has no nonidentity fixed points. A more general but equally checkable criterion involves the \textbf{secondary divisor}

\begin{equation}\label{Equation:SecondaryDivisor}
\gamma = \gcd(n,k-2j, l-2k+j, k-2l, j+l).
\end{equation}
The terms $k-2j, l-2k+j, -2l+k, -j-l$ are the differences of consecutive exponents on the coefficient $a$ in the cyclic relator $W = xa^jxa^{k-j}xa^{l-k}xa^{-l}$ for the shift extension. The dependence relation

$$
(k-2j) + (l-2k+j) + (k-2l) + (j+l) = 0.
$$
implies that any one of the four bracketed terms can be deleted when calculating $\gamma$. Thus $\gamma = \gcd(n,k-2j, k-2l, j+l)$. 

Since the shift $\theta_{G}$ arises from conjugation by $a$ in the shift extension $E = G \rtimes_{\theta_G} \Z_n \cong \pres{a,x}{a^n,W}$, the fixed point set for any power $\theta_G^d$ of the shift is expressible in terms of a centralizer in $E$:

$$
\mathrm{Fix}(\theta_G^d) = G \cap \mathrm{Cent}_E(a^d).
$$
In the proof of Theorem \ref{Theorem:secondary} below, we show that if $\gamma \neq 1$ then the shift extension $E$ admits an amalgamated free product decomposition with $\sgp{a} \cong \Z_n$ as a vertex group. We then apply the following general result, which is derived from \cite[Theorem 4.5]{MKS}. A proof using Bass-Serre theory is included.

\begin{Lemma}[Centralizer Lemma] \label{lemma:cent} Suppose that $\Pi = A \ast_C B$ is an amalgamated free product. If $a \in A$ and $g \in \mathrm{Cent}_\Pi(a)$, then either $g \in A$ or there exists $\alpha \in A$ such that $a \in \alpha C \alpha^{-1}$. 
\end{Lemma}

\proof The group $\Pi$ acts on the standard graph $T$ of $\Pi = A \ast_C B$ with orbit graph consisting of a single edge joining two distinct vertices. The vertices of $T$ are the cosets $wA$ and $wB$ where $w \in \Pi$. The edges of $T$ are the cosets $wC$ where $w \in \Pi$. The endpoints of the edge $wC$ are $wA$ and $wB$ and the action of $\Pi$ is the natural left action on cosets. Since $ga=ag$, the element $a \in A$ fixes both $1A$ and $gA$. If $g \not \in A$ then the fact that $T$ is a tree (see e.~g.~\cite[Theorem I.7.6]{DicksDunwoody}) implies that $a$ fixes the nontrivial geodesic joining the distinct vertices $1A$ and $gA$ in $T$, and hence $a$ fixes an edge adjacent to $1A$. This means that there exists $\alpha \in A$ such that $a\alpha C = \alpha C$ and so $a \in \alpha C \alpha^{-1}$. \qed

\begin{Theorem}[Secondary Divisor Criterion]\label{Theorem:secondary} Given $n,j,k,l$, let $G = G_n(x_0x_jx_kx_l)$ and let $\gamma =  \gcd(n,k-2j,l-2k+j, k-2l,j+l)$ be the secondary divisor. \begin{enumerate}
\item[(a)] If $\mathrm{Fix}(\theta_G^d) \neq 1$, then $\gamma \mid d$. 
\item[(b)] If $\gamma > 1$, then $G$ is infinite and $\mathrm{Fix}(\theta_G) = 1$.
\end{enumerate}
\end{Theorem}

\proof The substitution $u = xa^j$ reveals an amalgamated free product decomposition for the shift extension of $G$,
\begin{eqnarray*}
E &=& G \rtimes_\theta \Z_n = \pres{a,x}{a^n, xa^jxa^{k-j}xa^{l-k}xa^{-l}}\\
& \cong &\pres{a,u}{a^n, u^2a^{k-2j}ua^{l-k-j}ua^{-l-j}}\\
&\cong& \pres{a,b,u}{b=a^\gamma, a^n, u^2a^{k-2j}ua^{-(j+l)-(k-2l)}ua^{-(j+l)}}\\
&\cong& \Z_n \ast_{\Z_{n/\gamma}} D,
\end{eqnarray*}
where $\gamma = \gcd(n,k-2j, k-2l, j+l)$ and 

$$
D = \pres{b,u}{b^{n/\gamma}, u^2b^{(k-2j)/\gamma}ub^{(l-k-j)/\gamma}ub^{-(j+l)/\gamma}}.
$$
We have $D/\langle\langle b\rangle \rangle \cong \Z_4$, so the amalgamating subgroup $\sgp{b}= \sgp{a^\gamma} \cong \Z_{n/\gamma}$ is a proper subgroup of $D$. Now suppose that $\mathrm{Fix}(\theta_G^d) \neq 1$. Then we have $G \cap \mathrm{Cent}_E(a^d) \neq 1$ so using the fact that $G \cap \sgp{a} = 1$, the Centralizer Lemma implies that $a^d \in \sgp{a^\gamma}$. Since $\gamma \mid n$, this implies that $\gamma \mid d$. If we now suppose that $\gamma \neq 1$, then $\mathrm{Fix}(\theta_G) = 1$, as above, and moreover $\Z_{n/\gamma} = \sgp{a^\gamma}$ is a proper subgroup of $\Z_n = \sgp{a}$, so $E$ is infinite, whence the index $n$ subgroup $G$ is infinite. \qed

\begin{Corollary}[Primary Divisor Criterion]\label{Corollary:Primary} If $\gcd(n,j,k,l) \neq 1$, then the cyclically presented group $G = G_n(x_0x_jx_kx_l)$ is infinite and $\mathrm{Fix}(\theta_G) = 1$. Moreover if $\mathrm{Fix}(\theta_G^d) \neq 1$, then $\gcd(n,j,k,l) \mid d$.
\end{Corollary}

\proof One simply notes that any common divisor of $n,j,k,l$ is also a divisor of $\gamma$. \qed

\begin{Example}\label{Example:IsolatedExamples}\em For the following tuples $(n,j,k,l)$, the primary divisor $\gcd(n,j,k,l) = 1$ gives no information, but the secondary divisor $\gamma = \gcd(n,k-2j,k-2l,j+l)$ is greater than one, so the conclusions of Theorem \ref{Theorem:secondary} apply. Compare Tables \ref{Table:I} and \ref{Table:U}.
$$
\begin{array}{llll}
(10,3,6,1) & (12,1,2,9) & (16,3,6,1)&(20,3,6,1)\\
(24,1,2,15) & (24,3,6,1) & (24,1,2,19) &
\end{array}
$$
It is not difficult to show that if the primary divisor $\gcd(n,j,k,l) = 1$, then the secondary divisor $\gamma = \gcd(n,k-2j, k-2l, j+l)$ is equal to $1,2$, or $4$. We omit the details.
\end{Example}

\section{Conditions (A), (B), (C)}\label{Section:ABC}

Consider a cyclic presentation $\mathcal{P}_n(x_0x_jx_kx_l)$ for the group $G = G_n(x_0x_jx_kx_l)$ as in (\ref{Equation:P}). We now use the conditions (A), (B), (C) from Table \ref{Table:ABC} to simplify presentations for the shift extension $E = G \rtimes_{\theta_G} \Z_n$. We will use the following elementary fact without reference.

\begin{Lemma} Given $(n,j,k,l)$, if (A) and (C) are true, then both conditions in (A) are true and both conditions in (C) are true.
\end{Lemma}

We first consider the situation where condition (C) is true. 

\begin{Lemma}\label{Lemma:Cdecomp} Given $(n,j,k,l)$, let $G = G_n(x_0x_jx_kx_l)$ and assume that (C) is true. There exists an integer $p$ such that the shift extension $E = G \rtimes_{\theta_G} \Z_n$ admits a presentation of the form 
$$
E \cong \pres{a,z}{a^n, z^2a^{k-2p}z^2a^{-k-2p}} \cong \pres{a,u,z}{a^n, ua^kua^{-k}, ua^{2p} = z^2}
$$
where $G = \ker \nu$ for a retraction $\nu:E \ra \Z_n = \sgp{a}$ satisfying $\nu(a) = a$, $\nu(u) = 1$, and $\nu(z) = a^p$. If $l \equiv j+k$ mod $n$ in (C), then we take $p = j$. If $l \equiv j-k$ mod $n$ in (C), we take $p = -l$. %Moreover the secondary divisor $\gamma = \gcd(n,k+2p,k-2p)$.
\end{Lemma}

\proof If $l \equiv j+k$ mod $n$ then we have 
%\begin{eqnarray*}
%E &\cong& \pres{a,x}{a^n,xa^jxa^{k-j}xa^{l-k}xa^{-l}}\\
%&\cong& \pres{a,x}{a^n, xa^jxa^{k-j}xa^jxa^{-k-j}}\\
%&\cong& \pres{a,u,x}{a^n, u=xa^jxa^{-j}, ua^kua^{-k}}\\
%&\cong& \pres{a,u,z,x}{a^n, ua^kua^{-k}, z=xa^j, u=z^2a^{-2j}}\\
%&\cong& \pres{a,u,z}{a^n, ua^kua^{-k}, ua^{2j}=z^2}
%\end{eqnarray*}
\begin{eqnarray*}
E &\cong& \pres{a,x}{a^n,xa^jxa^{k-j}xa^{l-k}xa^{-l}}\\
&\cong& \pres{a,x}{a^n, xa^jxa^{k-j}xa^jxa^{-k-j}}\\
&\cong& \pres{a,x,z}{a^n,z=xa^j, z^2a^{k-2j}z^2a^{-k-2j}}\\
&\cong& \pres{a,z,u}{a^n, u=z^2a^{-2j}, ua^kua^{-k}}
\end{eqnarray*}
and so we take $p = j$ in this case. If $l \equiv j-k$, then a similar approach where $z=xa^{-l}$ and $u = z^2a^{2l}$ yields the desired result with $p=-l$. 
%\begin{eqnarray*}
%E &\cong& \pres{a,x}{a^n,xa^jxa^{k-j}xa^{l-k}xa^{-l}}\\
%&\cong& \pres{a,x}{a^n, xa^jxa^{k-j}xa^{j-2k}xa^{k-j}}\\
%&\cong& \pres{a,x}{a^n, xa^{k-j}xa^{j-2k}xa^{k-j}xa^j}\\
%&\cong& \pres{a,u,x}{a^n, u=xa^{k-j}xa^{j-k}, ua^{-k}ua^{k}}\\
%&\cong& \pres{a,u,z,x}{a^n, ua^kua^{-k}, z=xa^{k-j}, u=z^2a^{2(j-k)}}\\
%&\cong& \pres{a,u,z}{a^n, ua^kua^{-k}, ua^{2(k-j)}=z^2}\\
%&\cong& \pres{a,u,z}{a^n, ua^kua^{-k}, ua^{-2l}=z^2}
%\end{eqnarray*}
%so here $p=-l$. 
Since $\nu(a) = a$ and $\nu(x) = 1$, we have $\nu(u) = 1$ and $\nu(z)= a^p$ in both cases. \qed

%\begin{Lemma}\label{Lemma:C} Given $(n,j,k,l)$, let $G = G_n(x_0x_jx_kx_l)$ and assume that (C) is true. With a substitution of the form $u = xa^p$, the standard presentation for the split extension $E = G \rtimes_{\theta_G} \Z_n = \pres{a,x}{a^n, xa^jxa^{k-j}xa^{l-k}xa^{-l}}$ becomes
%
%$$
%E = G \rtimes_{\theta_G} \Z_n \cong \pres{a,u}{a^n, u^2 \alpha u^2 \beta}
%$$
%where $\alpha, \beta \in \sgp{a} \cong \Z_n$ are such that $\alpha\beta^{-1} = a^{\pm 2k}$. In this transformed presentation, 
%\begin{enumerate}
%\item[(a)] $\alpha = \beta \iff$ (A) is true, and
%\item[(b)] $\alpha$ or $\beta = 1 \iff$ (B) is true.
%\end{enumerate}
%\end{Lemma}
%
%\proof The condition (C) provides that $\l \equiv j \pm k$ mod $n$. If $l \equiv j+k$ mod $n$, then letting $u = xa^j$ we have 
%
%\begin{eqnarray*}
%E = G \rtimes_{\theta_G} \Z_n &=& \pres{a,x}{a^n, xa^jxa^{k-j}xa^{l-k}xa^{-l}}\\
%&\cong & \pres{a,x}{a^n, xa^jxa^{k-j}xa^{j}xa^{-j-k}}\\
%&\cong& \pres{a,u}{a^n, u^2a^{k-2j}u^2a^{-k-2j}}
%\end{eqnarray*}
%and so $\alpha = a^{k-2j}$ and $\beta = a^{-k-2j}$. It is now routine to verify that $\alpha = \beta \iff$ (A) is true and that $\alpha$ or $\beta = 1 \iff$ (B) is true. In the case where $l \equiv j-k$ mod $n$ as in (C), then the substitution $u = xa^{k-j}$ leads to the transformed relation $u^2\alpha u^2\beta$ where $\alpha = a^{2j-k}$ and $\beta = a^{2j-3k}$. \qed

\bigskip

The four cases in condition (B) can be reduced to a single case as follows.

\begin{Lemma}\label{Lemma:k=2j} If the parameters $(n,j,k,l)$ satisfy the condition (B), then there exists $c \in \Gamma_n$ such that $c(x_0x_jx_kx_l) = x_0x_{j'}x_{k'}x_{l'}$ where 

\begin{enumerate}
\item[(a)] $k' \equiv 2j'$ mod $n$;
\item[(b)] $\gcd(n,j',k',l') = \gcd(n,j,k,l)$;
\item[(c)] $\gcd(n,k'-2j',l'-2k'+j',k'-2l',j'+l') = \gcd(n,k-2j,l-2k+j,k-2l, j+l)$;
\item[(d)] The parameters $(n,j,k,l)$ satisfy (A) if and only if $(n,j',k',l')$ satisfy (A); and
\item[(e)] The parameters $(n,j,k,l)$ satisfy (C) if and only if $(n,j',k',l')$ satisfy (C).
\end{enumerate}
%$k' \equiv 2j'$ mod $n$. Moreover the primary divisors 
%
%$$
%\gcd(n,j,k,l)\ \ \mbox{and}\ \ \gcd(n,j',k',l')
%$$
%are the same and the secondary divisors
%
%$$
%\gamma = \gcd(n,k-2j,l-2k+j,k-2l, j+l)\ \ \mbox{and}\ \ \gamma' = \gcd(n,k'-2j',l'-2k'+j',k'-2l',j'+l')
%$$
%are the same.
\end{Lemma}

\proof All congruences are interpreted modulo $n$. If the parameters $(n,j,k,l)$ satisfy (B), then $k \equiv 2j$ or $k \equiv 2l$ or $j+l \equiv 0$ or $j+l \equiv 2k$. If $k \equiv 2j$ then there is nothing to prove and we take $c = 1 \in \Gamma_n$. 

If $k \equiv 2l$, then we take $c = \sigma$, so that $x_0x_{j'}x_{k'}x_{l'} = \sigma(x_0x_jx_kx_l) = x_0x_lx_kx_j$ so by setting $(n,j',k',l') = (n,l,k,j)$ we have $k' = k \equiv 2l = 2j'$. In addition $\gcd(n,j',k',l') = \gcd(n,l,j,k)$ and $$\gamma' = \gcd(n,k'-2j',l'-2k'+j',k'-2l',j'+l') = \gcd(n,k-2l,j-2k+l,k-2j,l+j) = \gamma.$$
One quickly checks that $(n,j',k',l') = (n,l,k,j)$ satisfies (A) if and only if $(n,j,k,l)$ satisfies (A), and similarly for (C).

If $j+l \equiv 0$, then we apply elements of $\Gamma_n$ as follows:
\begin{eqnarray*}
x_0x_jx_kx_l &=& x_0x_jx_kx_{-j}\\
&\stackrel{\theta_F^j}{\ra}& x_jx_{2j}x_{j+k}x_0\\
&\stackrel{\tau}{\ra}& x_0x_jx_{2j}x_{j+k}.
\end{eqnarray*}
Thus by setting $(n,j',k',l') = (n,j,2j,j+k)$ we obtain $k' = 2j = 2j'$. In addition, using $j+l \equiv 0$ we have $\gcd(n,j',k',l') = \gcd(n,j,2j,j+k) = \gcd(n,j,k) = \gcd(n,j,k,l)$ and 
\begin{eqnarray*}
\gamma' &=& \gcd(n,k'-2j',l'-2k'+j',k'-2l',j'+l')\\
&=& \gcd(n,0,j+k-4j+j,2j-2(j+k),2j+k)\\
&=& \gcd(n,k-2j,2k,2j+k)
\end{eqnarray*}
so that
$$
\gamma = \gcd(n,k-2j,l-2k+j,k-2l,j+l) = \gcd(n,k-2j,2k,k+2j,0) = \gamma'.
$$
Again using $j+l \equiv 0$, we note that
\begin{eqnarray*}
2k' \equiv 0 &\iff& 4j \equiv 0 \iff 2j-2l \equiv 0 \iff 2j \equiv 2l\\
2j' \equiv 2l' &\iff & 2j \equiv 2j +2k \iff 2k \equiv 0\\
l' \equiv j'+k' &\iff & j+k \equiv 3j \iff k \equiv 2j \iff l \equiv j-k\\
l' \equiv j'-k' &\iff & j+k \equiv -j \iff l \equiv j+k.
\end{eqnarray*}
Thus $(n,j',k',l')$ satisfies (A) (resp. (C)) if and only if $(n,j,k,l)$ satisfies (A) (resp. (C)).

And finally, if $j+l \equiv 2k$, then we apply elements of $\Gamma_n$ as follows:
\begin{eqnarray*}
x_0x_jx_kx_l &\stackrel{\tau^2}{\ra}& x_kx_lx_0x_j\\
&\stackrel{\sigma}{\ra}& x_kx_jx_0x_l\\
&\stackrel{\theta_F^{-k}}{\ra}& x_0x_{j-k}x_{-k}x_{l-k}
\end{eqnarray*}
Thus by setting $(n,j',k',l') = (n,j-k,-k,l-k)$, and by way of reducing to the previous case, we have $j'+l' = j-k+l-k \equiv j+l-2k \equiv 0$. In addition, we have $\gcd(n,j',k',l') = \gcd(n,j-k,-k,l-k) = \gcd(n,j,k,l)$ and 
\begin{eqnarray*}
\gamma' &=& \gcd(n,k'-2j',l'-2k'+j',k'-2l',j'+l')\\
&=& \gcd(n,-k-2(j-k),l-k+2k+j-k,-k-2(l-k),j-k+l-k)\\
&=& \gcd(n,,k-2j,j+l,k-2l,0) = \gamma.
\end{eqnarray*}
Again we note that
\begin{eqnarray*}
2k' \equiv 0 &\iff& 2k \equiv 0\\
2j' \equiv 2l' &\iff & 2j-2k \equiv 2l-2k \iff 2j \equiv 2l\\
l' \equiv j'+k' &\iff & l-k \equiv j-k-k \iff l \equiv j-k\\
l' \equiv j'-k' &\iff & l-k \equiv j-k+k \iff l \equiv j+k.
\end{eqnarray*}
Thus $(n,j',k',l')$ satisfies (A) (resp.~(C)) if and only if $(n,j,k,l)$ satisfies (A) (resp.~(C)). Thus the case $j+l \equiv 2k$ reduces to the case $j+k \equiv 0$ and the proof is complete. \qed

\begin{Lemma}\label{Lemma:B} Given $(n,j,k,l)$, let $G = G_n(x_0x_jx_kx_l)$ and assume that $k \equiv 2j$ mod $n$ as in (B). Let $\gamma = \gcd(n,k-2j,l-2k+j,k-2l,j+l)$ be the secondary divisor. With the substitution $u = xa^j$, the standard presentation for the shift extension $E = G \rtimes_{\theta_G} \Z_n$ takes the form $\pres{a,u}{a^n, u^3\alpha u\beta}$ where $\alpha = a^{l-3j}$ and $\beta = a^{-l-j}$ in $\sgp{a} \cong \Z_n$. Considering the parameters $(n,j,k,l)$ we have: 
\begin{enumerate}
\item[(a)] The condition (A) is true if and only if $\alpha = \beta^{\pm 1}$; 
\item[(b)] If (A) is true, then the element $\alpha \in \sgp{a} \cong \Z_n$ has order $n/\gamma$;% and it follows that if (C) is false, then $\alpha \neq 1$.
%in fact
%\begin{itemize}
%\item $\alpha = \beta^{-1} \iff$ (A1) and (B1) are true or (A2) and (B2) are true, and 
%\item $\alpha = \beta \iff$ (A1) and (B2) are true or (A2) and (B1) are true;
%\end{itemize}
\item[(c)] The condition (C) is true if and only if $\alpha = 1$ or $\beta = 1$;
%\item[(d)] If (C) is true, then $\gamma = \gcd(n,2k)$.
\end{enumerate}
\end{Lemma}

\proof Until otherwise stated, all congruences are to be considered modulo $n$. Using the substitution $u = xa^j$, we transform the standard presentation for the shift extension $E$ as follows. 
\begin{eqnarray*}
E &\cong & \pres{a,x}{a^n, xa^jxa^{k-j}xa^{l-k}xa^{-l}}\\
&\stackrel{k \equiv 2j}{\cong}& \pres{a,x}{a^n, xa^jxa^jxa^{l-2j}xa^{-l}}\\
&\stackrel{u=xa^j}{\cong}& \pres{a, u}{a^n, u^3a^{l-3j}ua^{-l-j}}.
\end{eqnarray*}
Thus $\alpha = a^{l-3j}$ and $\beta = a^{-l-j}$ and we verify (a) as follows:
\begin{eqnarray*}
\alpha = \beta^{-1} &\iff& l-3j \equiv j+l \iff 4j \equiv 0 \stackrel{k \equiv 2j}{\iff} 2k \equiv 0\\
\alpha = \beta &\iff& l-3j \equiv -l-j \iff 2l \equiv 2j.
\end{eqnarray*}

For statement (b), we claim that if either $2k \equiv 0$ or $2j \equiv 2l$, then the secondary divisor $\gamma$ has the following simplified description:
$$
\gamma = \gcd(n, k-2j, l-2k+j, j+l) = \gcd(n,j+l).
$$
This follows by noting that if $k \equiv 2j$ and $2k \equiv  0$ mod $n$, then working modulo any common divisor of $n$ and $j+l$, we have $l-2k+j \equiv 0$. And similarly, if $k \equiv 2j$ and $2j \equiv 2l$ mod $n$, then working modulo any common divisor of $n$ and $j+l$, we have $l-2k+j \equiv l-3j \equiv l -2l -j \equiv -l-j \equiv 0$. From this it follows that the order of $\alpha = \beta^{\pm 1}$ is $n/\gcd(n,j+l) = n/\gamma$. 

For statement (c), $\alpha = 1 \iff l \equiv 3j \iff l \equiv j+2j \equiv j+k$, as in (C). %and in this case
%$$
%\gamma = \gcd(n,k-2j,j-2k+l,k-2l,j+l) = \gcd(n,0,4j-2k,k-6j,j+l) = \gcd(n,2k).
%$$
Also $j-k \equiv j-2j \equiv -j$, so $\beta = 1 \iff j+l \equiv 0 \iff l \equiv -j \equiv j-k$, as in (C). %and in this case
%$$
%\gamma = \gcd(n,k-2j,j-2k+l,k-2l,j+l) = \gcd(n,0,-2k,4j,0) = \gcd(n,2k).
%$$
This completes the proof. \qed

\section{C(4)-T(4) Cases}\label{Section:SmallCancellation}

In this section we consider those cases that satisfy the C(4)-T(4) small cancellation conditions. We begin with the observation that with just one exception, finiteness and asphericity are mutually exclusive.

\begin{Lemma}\label{Lemma:FiniteAsph} Assume that $\mathcal{P} = \mathcal{P}_n(x_0x_jx_kx_l)$ is combinatorially aspherical and let $G = G_n(x_0x_jx_kx_l)$. 
\begin{enumerate}
\item[(a)] If $G$ has nontrivial torsion, then $k \equiv 0$ and $j\equiv l$ mod $n$, so (A) and (C) are true.
\item[(b)] If $G$ is finite, then $n=1$, so (A), (B), and (C) are true, $\mathcal{P} = \pres{x_0}{x_0^4}$, and $G \cong \Z_4$. 
\end{enumerate}
\end{Lemma}

\proof Since $\mathcal{P}$ is combinatorially aspherical and $G$ has torsion, \cite[Theorem 3]{Hueb} implies that the relator $w = x_0x_jx_kx_l$ or one of its shifts must be a proper power in the free group $F$ with basis $x_0,x_1, \ldots, x_{n-1}$. This implies that $k \equiv 0$ and $j \equiv l$ mod $n$ so that $G \cong G_n((x_0x_j)^2)$. This proves (a). If $G$ is finite, then since $G$ is nontrivial (having a $\Z_4$ homomorphic image), we know that $k \equiv 0$ mod $n$ and (C) is true so by Lemma \ref{Lemma:Cdecomp}, the shift extension has a presentation $E = G \rtimes_{\theta_G} \Z_n \cong \pres{a,u,z}{a^n, u^2, ua^{2p} = z^2}$ and so $E$ contains the free product $\pres{a,u}{a^n,u^2} \cong \Z_n \ast \Z_2$ as a subgroup. If $G$ is finite, then $E$ must be finite and so $n=1$, as in (b). \qed

\begin{Lemma}\label{Lemma:SmCanc} If (B) and (C) are false, then $\mathcal{P}_n(x_0x_jx_kx_l)$ satisfies the C(4)-T(4) small cancellation conditions and so is combinatorially aspherical.
\end{Lemma}

\proof If (B) and (C) are false, then $j, k-j,l-k,-l$ are pairwise distinct modulo $n$. It follows easily from Lemma \ref{Lemma:SmCancCrit} that $\mathcal{P}_n(x_0x_jx_kx_l)$ satisfies the C(4)-T(4) small cancellation conditions. Alternatively, \cite[Theorem 2]{BBP} provides that the relative presentation $\pres{\Z_n,x}{xa^jxa^{k-j}xa^{l-k}xa^{-l}}$ for the shift extension is aspherical and so $\mathcal{P}_n(x_0x_jx_kx_l)$ is combinatorially aspherical by Theorem \ref{Theorem:AsphTransf}. \qed

\begin{Theorem}\label{Lemma:Decomposable} Given $(n,j,k,l)$, if (C) is true, then the following are equivalent:
\begin{enumerate}
\item[(i)] (A) is true,
\item[(ii)] $\mathcal{P}_n(x_0x_jx_kx_l)$ satisfies the C(4)-T(4) small cancellation conditions,
\item[(iii)] $\mathcal{P}_n(x_0x_jx_kx_l)$ is combinatorially aspherical,
\item[(iv)] $\Z_n$ acts freely via the shift on the nonidentity elements of $G_n(x_0x_jx_kx_l)$.
\end{enumerate}
\end{Theorem}

\proof Assuming that (C) is true, the word $w = x_0x_jx_kx_l$ takes the form $w = x_0x_jx_kx_{j \pm k}$. We first show that (i) implies (ii), so suppose that (A) is true. Then $2k \equiv 0$ mod $n$ and so $w = x_0x_jx_kx_{j+k}$. Up to cyclic shifts, the only length two cyclic subwords of $w$ are $x_0x_j$ and $x_jx_k$, so by Lemma \ref{Lemma:SmCancCrit} it suffices to show that neither of these words is a piece. So suppose that the initial subword $x_0x_j$ of $w = x_0x_jx_kx_{j+k}$ is equal to an initial subword of a cyclic permutation $v$ of a shift $\theta_F^i(w) = x_ix_{i+j}x_{i+k}x_{i+j+k}$ of $w$. We are to show that $w = v$. There are four cases to check. If $v = \theta_F^i(w) = x_ix_{i+j}x_{i+k}x_{i+j+k}$, then $i=0$ so $v = x_0x_jx_kx_{j+k} = w$. If $v = x_{i+j}x_{i+k}x_{i+j+k}x_i$, then $x_{i+j}x_{i+k} = x_0x_j$, so $i+j \equiv 0$ and $i+k \equiv j$ (congruences modulo $n$). Then $i+j+k \equiv k$ and the fact that $2k \equiv 0$ mod $n$ implies that $i \equiv i+2k \equiv j+k$, so $w = v$. Similar arguments apply to the remaining two cyclic permutations of $\theta_F^i(w)$ to show that $x_0x_j$ is not a piece. One argues similarly to show that $x_jx_k$ is not a piece by comparing the cyclic permutation $x_jx_kx_{j+k}x_0$ of $w$ with the cyclic permutations $v$ of the shift $\theta_F^i(w) = x_ix_{i+j}x_{i+k}x_{i+j+k}$. It follows that no length two cyclic subword of $w$ is a piece and so $\mathcal{P}_n(x_0x_jx_kx_l)$ satisfies C(4)-T(4) by Lemma \ref{Lemma:SmCancCrit}. 

As previously noted, C(4)-T(4) presentations are combinatorially aspherical, so (ii) implies (iii), and Theorem \ref{Theorem:Free} provides that (iii) implies (iv), so it remains to show that if (A) is false and (C) is true, then some power of the shift $\theta_G$ on $G = G_n(x_0x_jx_kx_l)$ has a nonidentity fixed point. In the presentation for the shift extension $E = G \rtimes_{\theta_G} \Z_n$ from Lemma \ref{Lemma:Cdecomp} we have $a,u \in E$ satisfying $ua^kua^{-k} = 1$ where $u$ does not lie in the normal closure $\<\<a\>\>_E$ of $a$ because $u \neq 1$ in $E/\<\<a\>\>_E \cong \pres{z,u}{u^2,u = z^2} \cong \Z_4$. Thus it follows that the element $g = u\nu(u)^{-1}$ is a nonidentity element of $G = \ker \nu$ where $\nu: E \ra \Z_n$ is a retraction satisfying $\nu(a) = a$. Further, the fact that $ua^kua^{-k} = 1$ implies that $u \in \mathrm{Cent}_E(a^{2k})$ and hence $1 \neq g \in G \cap \mathrm{Cent}_E(a^{2k}) = \mathrm{Fix}(\theta_G^{2k})$. Since (A) is false, we also have $2k \not \equiv 0$ mod $n$ and so the group $\sgp{a} \cong \Z_n$ does not act freely via the shift. \qed

%\begin{Theorem}\label{Theorem:FiniteACnotBC} Given $(n,j,k,l)$, let $G = G_n(x_0x_jx_kx_l)$. If (A) and (C) are true or (B) and (C) are false, then the following are equivalent:
%\begin{enumerate}
%\item[(i)] $G$ is finite;
%\item[(ii)] $n=1$;
%\item[(iii)] $\theta_G = 1$;
%\item[(iv)] $\mathrm{Fix}(\theta_G) \neq 1$.
%\end{enumerate} 
%\end{Theorem}

%\proof These are the C(4)-T(4) cases, so the result follows from Lemma \ref{Lemma:FiniteAsph}. \qed

\section{FxT}\label{Section:FxT}

In this section we consider the case where (C) is true and (A) is false. Note that $\mathcal{P}_n(x_0x_jx_kx_l)$ is not combinatorially aspherical and the shift action is not free by Theorem \ref{Lemma:Decomposable}.

\begin{Theorem}\label{Theorem:CFix} Given $(n,j,k,l)$, let $G = G_n(x_0x_jx_kx_l)$. If (C) is true and (A) is false, then the following are equivalent:
\begin{enumerate}
\item[(i)] The shift $\theta_G$ has a nonidentity fixed point;
\item[(ii)] The group $D = \pres{a,u}{a^n,ua^kua^{-k}}$ is finite;
\item[(iii)] $\gcd(n,2k) = 1$.
\end{enumerate} 
\end{Theorem}

\proof By Lemma \ref{Lemma:Cdecomp} the shift extension $E = G \rtimes_{\theta_G} \Z_n$ has a presentation 
$$E \cong \pres{a,u,z}{a^n, ua^kua^{-k}, z^2 = ua^{2p}} \cong \pres{D,z}{z^2=ua^{2p}}.$$
We have the retraction $\nu: E \ra \Z_n = \sgp{a}$ with kernel $G = \ker \nu$ where $\nu(a) = a$, $\nu(u) = 1$, and $\nu(z) = a^p$. The restriction $\nu|_D: D \ra \Z_n$ has cyclically presented kernel $H = \ker \nu \cap D \cong G_n(u_0u_k)$. (See \cite[Theorem 2.3]{BShift}.) Letting $\kappa = \gcd(n,k)$, it is routine to verify that 
$$
H \cong \ast_{i=1}^\kappa G_n(u_0u_1) \cong \left\{
\begin{array}{rl}
\ast_{i=1}^\kappa \Z & \mbox{if $n$ is even}\\
\ast_{i=1}^\kappa \Z_2 & \mbox{if $n$ is odd}
\end{array} \right.
$$
Since $D$ is finite if and only if $H$ if finite, this establishes (ii) $\iff$ (iii). 

Next we show that (i) $\Rightarrow$ (iii). The split extension $E = G \rtimes_{\theta_G} \Z_n$ decomposes as $E = D \ast_{\sgp{ua^{2p} = z^2}} \sgp{z} = D \ast_Y \sgp{z}$ where $Y = \sgp{ua^{2p}} = \sgp{z^2} = D \cap \sgp{z}$. If we assume that $1 \neq g \in \mathrm{Fix}(\theta_G) = G \cap \mathrm{Cent}_E(a)$, as in (i), then since $a \in D$, the Centralizer Lemma \ref{lemma:cent} implies that either $g \in D$ or there exists $\alpha \in D$ such that $a \in \alpha Y \alpha^{-1}$. Either way, we claim that $n$ is odd and $\gcd(n,k) = 1$. If $g \in D$, then $1 \neq g \in \mathrm{Fix}(\theta_H)$ and so $\gcd(n,k) = 1$ by an obvious length two analogue of Corollary \ref{Corollary:Primary}. Further, $n$ must be odd because otherwise $H \cong G_{n}(u_0u_1) \cong \Z$ and $\theta_H = -1$ is fixed point free. So suppose that $a \in \alpha Y \alpha^{-1}$ where $\alpha \in D$. In this case $a$ lies in the normal closure $\<\<Y\>\>_D$ of $Y$ in $D$ and so the fact that
$$
D/\<\<Y\>\>_D \cong \pres{a,u}{a^n,ua^kua^{-k}, u=a^{-2p}} \cong \pres{a}{a^n,a^{-4p}} \cong \pres{a}{a^{\gcd(n,4p)}}
$$
implies that $\gcd(n,4p) = 1$, whence $n$ is odd. Working in $D$, the element $a$ is a power of $\alpha ua^{2p}\alpha^{-1}$ where $\alpha \in D$, so that $\alpha ua^{2p}\alpha^{-1} \in \mathrm{Cent}_D(a)$. Letting $b = a^k$, we have
$$
D = \pres{a,u}{a^n,ua^kua^{-k}} \cong \pres{a}{a^n} \ast_{a^k = b} \pres{b,u}{b^{n/\gcd(n,k)}, ubub^{-1}} = \Z_n \ast_B C
$$
where $C = \pres{b,u}{b^{n/\gcd(n,k)},ubub^{-1}}$ and $B = \sgp{b} = \sgp{a^k} \cong \Z_{n/\gcd(n,k)}$. By the Centralizer Lemma \ref{lemma:cent} it follows that either $\alpha ua^{2p}\alpha^{-1} \in \sgp{a}$ or else $a \in \beta B \beta^{-1} = B = \sgp{a^k}$ for some $\beta \in \sgp{a}$. Now the first of these conditions cannot be true because $D/\<\<a\>\>_D \cong \pres{u}{u^2}$ so $u \not \in \<\<a\>\>_D$ and hence $ua^{2p} \not \in \<\<a\>\>_D$. Thus it follows that $a \in \sgp{a^k}$ and so $\gcd(n,k) = 1$. 

It remains to show that if $n$ is odd and $\gcd(n,k) = 1$, then $\theta_G$ has a nonidentity fixed point. For this we have $1 \neq u \in \ker \nu = G$ and $ua^k ua^{-k} = 1$ and hence $a^{2k}ua^{-2k} =u$. Since $\gcd(n,2k) = 1$, it follows that $u \in G \cap \mathrm{Cent}_E(a^{2k}) = G \cap \mathrm{Cent}_E(a) = \mathrm{Fix}(\theta_G)$. In fact, $D \cong \Z_2 \times \Z_n \cong \Z_{2n}$. \qed

\begin{Theorem}\label{Theorem:FiniteCnotA} Given $(n,j,k,l)$, let $G = G_n(x_0x_jx_kx_l)$. If (C) is true and (A) is false, then the following are equivalent:
\begin{enumerate}
\item[(i)] $E = G \rtimes_{\theta_G} \Z_n \cong \Z_{4n}$;
\item[(ii)] $G \cong \Z_4$;
\item[(iii)] $G$ is finite;
\item[(iv)] $\gcd(n,2k) =1$ and either
\begin{itemize}
\item $l \equiv j+k$ mod $n$ and $\gcd(n,j) = 1$, or
\item $l \equiv j-k$ mod $n$ and $\gcd(n,l) = 1$.
\end{itemize}
\end{enumerate}
\end{Theorem}

\proof The implications (i) $\Rightarrow$ (ii) $\Rightarrow$ (iii) being obvious, assume that $G$ is finite, so that by Lemma \ref{Lemma:Cdecomp} the shift extension $$E = G \rtimes_{\theta_G} \Z_n \cong \pres{a,u,z}{a^n,ua^kua^{-k},z^2=ua^{2p}} \cong D \ast_{ua^{2p}=z^2} \sgp{z}$$ is finite where $D \cong \pres{a,u}{a^n,ua^kua^{-k}}$. Thus $D$ is finite and so $\gcd(n,2k) = 1$ by Theorem \ref{Theorem:CFix}, as in (iv). Further,
\begin{eqnarray*}
E &\cong& \pres{a,z}{a^n,z^2a^{k-2p}z^2a^{-k-2p}}\\
&\cong & \pres{a,v,z}{a^n, v^2a^{-2k}, v = z^2a^{k-2p}}\\
&\cong& \pres{a,v}{a^n,v^2a^{-2k}} \ast_{va^{2p-k} = z^2} \sgp{z}.
\end{eqnarray*}
Since $\gcd(n,2k) = 1$, it follows that $v$ is a generator of $\pres{a,v}{a^n, v^2a^{-2k}} = \sgp{v} \cong \Z_{2n}$ and so $E \cong \Z_{2n} \ast_{va^{2p-k}=z^2}\sgp{z}$. Now $G$ is finite if and only if the element  $va^{2p-k} \in \sgp{v} \cong \Z_{2n}$ has order $2n$, in which case $E \cong \Z_{4n}$. It remains to verify that $G$ is finite if and only if one of the itemized conditions in (iv) is true. 

Let $k'$ be a multiplicative inverse for $2k$ modulo $n$ so that $2kk' \equiv 1$ mod $n$. Then $v^{2k'} = a^{2kk'} = a$ and hence the order of the element
$$
va^{2p-k} = v^{1+2k'\cdot 2p - 2kk'} \in \sgp{v} \cong \Z_{2n}
$$
is $2n/g$ where $$g = \gcd(2n,1+4pk' - 2kk') = \gcd(n,1+4pk' - 2kk') = \gcd(n,1+4pk' - 1) = \gcd(n,p).$$ Thus $G$ is finite if and only if $\gcd(n,p) = 1$. Since (C) is true we have $l \equiv j \pm k$ mod $n$. If $l \equiv j+k$ then by Lemma \ref{Lemma:Cdecomp} we have $p=j$, whereas if $l \equiv j-k$ then $p = -l$. The result follows. \qed

\begin{Example}\label{Example:CFinite} \em  If (C) is true and (A) is false, then $G \cong \Z_4$ can occur regardless of whether (B) is true or not, as the following exemples demonstrate.
\begin{itemize}
\item If $G = G_n(x_0x_1x_2x_3)$ and $n > 4$, then (A) is false, (B) and (C) are true, and $G\ \mbox{is finite} \iff G \cong \Z_4 \iff n$ is odd.
\item If $G = G_n(x_0x_1x_3x_4)$ and $n > 6$, then (A) and (B) are false, (C) is true, and $G\ \mbox{is finite} \iff G \cong \Z_4 \iff \gcd(n,6) = 1$.
\end{itemize}

%If (A) and (C) are true, then $\mathcal{P}_n(x_0x_jx_kx_l)$ is combinatorially aspherical by Theorem \ref{Lemma:Decomposable} and so $G_n(x_0x_jx_kx_l)$ is finite if and only if $n=1$ by Lemma \ref{Lemma:FiniteAsph}.
\end{Example}

\section{TTF}\label{Section:TTF}

In this section we consider the case where (A) and (B) are true and (C) is false.

\begin{Theorem}\label{Theorem:MZ} Given $(n,j,k,l)$, assume that (A) and (B) are true and (C) is false.
\begin{enumerate}
%\item[(a)] $\mathcal{P}_n(x_0x_jx_kx_l)$ is not combinatorially aspherical.
\item[(a)] The shift action on $G = G_n(x_0x_jx_kx_l)$ is not free, so $\mathcal{P}_n(x_0x_jx_kx_l)$ is not combinatorially aspherical by Theorem \ref{Theorem:Free}.
\item[(b)] The following are equivalent:
\begin{enumerate}
\item[(i)] $G$ is finite;
\item[(ii)] The secondary divisor $\gamma = \gcd(n,k-2j,j-2k+l,k-2l,j+l) = 1$;
\item[(iii)] The shift $\theta_G$ has a nonidentity fixed point.
\end{enumerate}
\item[(c)] If $G$ is finite, then $G$ is solvable.
\end{enumerate}
\end{Theorem}

\proof Let $w = x_0x_jx_kx_l$. By Lemma \ref{Lemma:k=2j}, there exists $c \in \Gamma_n$ such that $c(w) = x_0x_{j'}x_{k'}x_{l'}$ where $k' \equiv 2j'$ mod $n$ and where the primary and secondary divisors, as well as the status of the conditions (A) and (C), are identical for $w$ and for $c(w)$. Further, the group structure, asphericity status, and shift dynamics are identical for the presentations $\mathcal{P}_n(w)$ and $\mathcal{P}_n(c(w))$ by Theorem \ref{Theorem:DUI}. Thus we may assume that the parameters $(n,j,k,l)$ are such that $k \equiv 2j$ mod $n$, (A) is true, and (C) is false.

By Lemma \ref{Lemma:B} the shift extension $E = G \rtimes_{\theta_G} \Z_n$ has a presentation $\pres{a,u}{a^n, u^3\alpha u\beta}$ where $1 \neq \alpha = \beta^{\pm 1}$ has order $n/\gamma$ in $\sgp{a} \cong \Z_n$, where $\gamma$ is is the secondary divisor. We verify the claims (a)-(c) in the cases $\alpha = \beta^{-1}$ and $\alpha = \beta$ separately. Recall that $G$ is the kernel of a retraction $\nu: E \ra \sgp{a} \cong \Z_n$. 

Suppose first that $1 \neq \alpha = \beta^{-1}$ so that we have an amalgamated free product decomposition
$$
E \cong \pres{a,u}{a^n, u^3\alpha u \alpha^{-1}} \cong \Z_n \ast_{\sgp{\alpha}} M
$$
where $M$ is the group with presentation $M \cong \pres{\alpha,u}{\alpha^{n/\gamma},u^3\alpha u \alpha^{-1}}$. By  \cite[Lemma 2.2]{BW2}, the group $M$ is the split metacyclic group $M \cong \sgp{u} \rtimes \sgp{\alpha}$, where $u$ has order $\mu = 3^{n/\gamma}-(-1)^{n/\gamma}$ and $\sgp{u} \cap \sgp{\alpha} = 1$ in $E$. The order $\mu$ is divisible by four and since (C) is false, Lemma \ref{Lemma:B}(c) provides that $n/\gamma > 1$, so $\mu > 4$, which means that the element $v = u^{\mu/4} \in M$ is a nontrivial element of $E$ that lies outside of $\sgp{\alpha} = \sgp{a} \cap M$. Now
$$
\alpha v \alpha^{-1} = \alpha u^{\mu/4} \alpha^{-1} = u^{-3\mu/4} = u^{(1-4)\mu/4} = v u^{-\mu} = v,
$$
and so $v \in \mathrm{Cent}_E(\alpha)$. Since $v \not \in \sgp{a}$ it follows that $g = v\nu(v)^{-1} \in G \cap \mathrm{Cent}_E(\alpha)$ is nontrivial. Since $\alpha \neq 1$ in $\sgp{a} \cong \Z_n$, we have shown that a nonidentity power of $a$ centralizes a nonidentity element of $G$ in the shift extension $E$. Thus the shift action on $G$ is not free, as in (a). Since $\sgp{\alpha}$ is a proper subgroup of $M$, we have that
$$
G\ \mbox{is finite} \iff E\ \mbox{is finite} \iff E=M \iff \sgp{a} = \sgp{\alpha} \iff \gamma = 1,
$$
in which case $G$ is a subgroup of the metacyclic group $M$ and so $G$ is metacyclic, hence solvable, as in (c). If $\gamma = 1$ then $1 \neq g \in G \cap \mathrm{Cent}_E(\alpha) =  G \cap \mathrm{Cent}_E(a) = \mathrm{Fix}(\theta_G)$, whereas if $\gamma \neq 1$ then $\mathrm{Fix}(\theta_G) = 1$ by Theorem \ref{Theorem:secondary}. This completes the proof of of the theorem in the case where $\beta = \alpha^{-1}$.

Assume that $1 \neq \alpha = \beta $ has order $n/\gamma$ in the amalgamated free product
$$
E \cong \pres{a,u}{a^n, u^3\alpha u \alpha} \cong \Z_n \ast_{\sgp{\alpha}} \Delta
$$
where $\Delta = \pres{\alpha,u}{\alpha^{n/\gamma}, u^3\alpha u\alpha}$. Working in $\Delta$ we have 
\begin{eqnarray*}
u^2\alpha \cdot u\alpha u &=& 1\ \mbox{and}\\
\alpha u^2 \cdot u \alpha u &=& 1,
\end{eqnarray*}
which implies that $u^2 \in \mathrm{Cent}_E(\alpha)$ and so $g = u^2\nu(u^2)^{-1} \in G \cap \mathrm{Cent}_E({\alpha})$.
Since $u^2 \not \in \<\<a\>\>_E$ we know that $1 \neq g \in G$. Again the fact that (C) is false means that $1 \neq \alpha \in \sgp{a} \cong \Z_n$, so we conclude that a nonidentity power of $a$ centralizes a nonidentity element of $G$. Thus the shift action on $G$ is not free, as in (a). 

The central quotient $\Delta/\sgp{u^2} = \pres{\alpha,u}{\alpha^{n/\gamma},u^2,u^3\alpha u \alpha} \cong \pres{\alpha,u}{\alpha^{n/\gamma},u^2, (u\alpha)^2}$ is the dihedral group $D_{n/\gamma}$ of order $2n/\gamma$ and we have the five-term homology sequence
$$
H_2\Delta \ra H_2D_{n/\gamma} \ra \sgp{u_2} \ra H_1\Delta \ra H_1D_{n/\gamma} \ra 0.
$$
The group $\Delta$ has a $2\times 2$ relation matrix with determinant $4n/\gamma \neq 0$ so $H_2\Delta = 0$ and $|H_1\Delta| = 4n/\gamma$. For the dihedral group we have $|H_1D_{n/\gamma}| - |H_2D_{n/\gamma}| = 2$ and so the five term sequence implies that
$$
0 = |H_2\Delta| - |H_2D_{n/\gamma}| +  |\sgp{u_2}| - |H_1\Delta| + |H_1D_{n/\gamma}| = |\sgp{u^2}| + 2-4n/\gamma.
$$
Thus $\Delta$ is a central extension of the dihedral group of order $2n/\gamma$ with cyclic kernel $\sgp{u^2}$ of order $-2+4n/\gamma$. In particular, $\Delta$ is solvable. As before, $\sgp{\alpha}$ is a proper subgroup of $\Delta$ and so 
$$
G\ \mbox{is finite} \iff E\ \mbox{is finite} \iff E=\Delta \iff \sgp{a} = \sgp{\alpha} \iff \gamma = 1,
$$
in which case $G$ is a subgroup of the solvable group $\Delta$ and so $G$ is solvable, as in (c). If $\gamma = 1$ then $1 \neq g \in G \cap \mathrm{Cent}_E(\alpha) =  G \cap \mathrm{Cent}_E(a) = \mathrm{Fix}(\theta_G)$, whereas if $\gamma \neq 1$ then $\mathrm{Fix}(\theta_G) = 1$ by Theorem \ref{Theorem:secondary}. This completes the proof of of the theorem in the case where $\beta = \alpha$. \qed

%
%\begin{enumerate}
%\item[\bf (M)] At least one of the following:
%\begin{enumerate}
%\item[\bf(M1)] $2k \equiv 0$ and either $k \equiv 2j$ or $k \equiv 2l$, or
%\item[\bf(M2)] $2l \equiv 2j$ and either $2k \equiv j+l$ or $j+l \equiv 0$
%\end{enumerate}
%\item[\bf (Z)]  At least one of the following:
%\begin{enumerate}
%\item[\bf(Z1)] $2j \equiv k \equiv 2l$, or
%\item[\bf(Z2)] $j+l \equiv 2k \equiv 0$.
%\end{enumerate}
%\end{enumerate}

\section{FTF}\label{Section:FTF}

Thus far we have considered seven of the eight possible combinations of the conditions (A), (B), (C). The remaining possibility, when (B) is true and both (A) and (C) are false, is the most complex and the most interesting. Here we will invoke the main result from \cite{BBP} to conclude that in most cases, the cyclic presentation $\mathcal{P}_n(x_0x_jx_kx_l)$ is combinatorially aspherical (see Theorem \ref{Theorem:Asphericity}), in which case we know all that we need to know. Fortunately, the remaining nonaspherical cases are just a few in number and so it is possible to analyze these individually, in some cases via computational means. Exemplars for \textbf{isolated} (I) and \textbf{unresolved} (U) conditions appear in Tables \ref{Table:I} and \ref{Table:U} of Section \ref{Section:Introduction}.
One checks case-by-case that if (I) or (U) is true for parameters $(n,j,k,l)$, then the primary divisor $\gcd(n,j,k,l) = 1$, (B) is true, and both (A) and (C) are false. Recall that each individual condition within (I) and (U) describes a complete orbit of a word $w \in \Phi_n$ under the action of $\Gamma_n = D_4 \times (\Z_n \rtimes \Z_n^\ast)$ as in Lemma \ref{Lemma:DU}. Thus it follows from Theorem \ref{Theorem:DUI} and Corollary \ref{Corollary:DUDyn} that each of these ten individual isolated or unresolved conditions refers to exactly one group up to isomorphism, asphericity status, and shift dynamics.

In order to deal with cases where the primary divisor $c = \gcd(n,j,k,l)$ is not equal to one, we introduced the generalized conditions (I$^\ast$) and (U$^\ast$). The asphericity status for $\mathcal{P}_n(x_0x_jx_kx_l)$ in types (I$^\ast$) and (U$^\ast$) is entirely determined by that of (I) and (U) because of the following well-known result. 

\begin{Lemma}[Free Product]\label{Lemma:FreeProduct} The cyclic presentation $\mathcal{P}_{nd}(x_0x_{jd}x_{kd}x_{ld})$ is combinatorially aspherical if and only if the cyclic presentation $\mathcal{P}_n(x_0x_jx_kx_l)$ is combinatorially aspherical.
\end{Lemma}

\proof As in \cite{EdjvetIrred}, $\mathcal{P}_{nd}(x_0x_{jd}x_{kd}x_{ld})$ decomposes as a disjoint union of $d$ subpresentations that can each be identified with $\mathcal{P}_n(x_0x_jx_kx_l)$, so this follows from \cite[Theorem 4.2]{CCH}. \qed

\bigskip

Given the cyclic presentation $\mathcal{P}_n(x_0x_jx_kx_l)$ for the group $G = G_n(x_0x_jx_kx_l)$, the shift extension admits the relative presentation 
$$
E = G \rtimes_{\theta_G} \Z_n \cong \pres{\Z_n,x}{xa^jxa^{k-j}xa^{l-k}xa^{-l}}
$$ 
with finite cyclic coefficient group $\sgp{a} \cong \Z_n$. This relative presentation is aspherical if and only if $\mathcal{P}_n(x_0x_jx_kx_l)$ is combinatorially aspherical by Theorem \ref{Theorem:AsphTransf}. A classification of asphericity for relative presentations of the form $\pres{H,x}{xh_1xh_2xh_3xh_4}$ where the coefficients $h_i$ are taken from an arbitrarily given group $H$ is given in \cite{BBP}, modulo a handful of \textbf{exceptional} cases where the asphericity status is unresolved. The analysis in \cite{BBP} first deals with the case where the coefficients $h_i \in H$ are pairwise distinct \cite[Theorem 2]{BBP}, in which case the relative presentation is aspherical. This is essentially a C(4)-T(4) small cancellation situation as in \cite{Lyndon66} and corresponds to Lemma \ref{Lemma:SmCanc} above. The situation where $h_1 = h_3$ or $h_2 = h_4$ is dealt with in \cite[Theorem 3]{BBP} and corresponds to Theorem \ref{Lemma:Decomposable} above. Thus we are left with the situation where two consecutive coefficients coincide in the relator, viewed as a cyclic word. Translated into the current context, where $H = \sgp{a} \cong \Z_n$ and the relator has the form $xa^jxa^{k-j}xa^{l-k}xa^{-l}$, this simply means that (B) is true. As in Lemma \ref{Lemma:B}, a linear substitution $u = xa^q$ now transforms the relative presentation $\mathcal{E} = \pres{\Z_n,x}{xa^jxa^{k-j}xa^{l-k}xa^{-l}}$ for the shift extension $E$ into one of the form $\mathcal{Q} = \pres{\Z_n,u}{u^3 \alpha u \beta}$ where $\alpha, \beta \in \sgp{a} \cong \Z_n$. Because of the simple nature of this substitution, the cellular model of the relative presentation $\mathcal{Q}$ is homotopy equivalent to that of $\mathcal{E}$ and so these relative presentations have the same asphericity status.

The statement of \cite[Theorem 4]{BBP} classifies asphericity of relative presentations $$\mathcal{Q} = \pres{H,u}{u^3 \alpha u \beta}$$ where $\alpha,\beta  \in H$ under the assumption that $2 \leq \mathrm{o}(\alpha) \leq \mathrm{o}(\beta )$ and where $\mathrm{o}(-)$ refers to the order of an element in $H$. Here we present a symmetrized version of the result that does not include this restriction on the orders of $\alpha$ and $\beta$ and which is updated in light of recent work. The following circumstances were termed \textbf{exceptional} in \cite{BBP}.

\begin{enumerate}
\item[\bf E1] $\alpha = \beta ^2$ or $\beta  = \alpha^2$ and $\sgp{\alpha,\beta } \cong \Z_6$
\item[\bf E2] $\alpha = \beta ^3$ or $\beta  = \alpha^3$ and $\sgp{\alpha,\beta } \cong \Z_6$
\item[\bf E3] $\alpha = \beta ^4$ or $\beta  = \alpha^4$ and $\sgp{\alpha,\beta } \cong \Z_6$
\item[\bf E4] $\{\mathrm{o}(\alpha), \mathrm{o}(\beta )\} = \{2,4\}$ and $\sgp{\alpha,\beta } \cong \Z_2 \times \Z_4$
\item[\bf E5] $\{\mathrm{o}(\alpha), \mathrm{o}(\beta )\} = \{2,5\}$ and $\sgp{\alpha,\beta } \cong \Z_2 \times \Z_5$
\end{enumerate}
As in other asphericity studies \cite{ECap,HM}, these exceptional cases were \textbf{unresolved} for asphericity status in \cite{BBP}. In recent work, Aldwaik and Edjvet \cite{AE14} have shown that $\mathcal{Q}$ is aspherical in the previously unresolved cases \textbf{E4} and \textbf{E5}. Also, Williams \cite{GW15} has informed us of calculations in GAP \cite{GAP} showing that the group $E = \pres{a,u}{a^6, u^3a^3ua}$ is finite of order $24\,530\,688 = 2^8 \cdot 3^4 \cdot 7 \cdot 13^2$ with second derived subgroup $E\"$ isomorphic to the simple group $\mathrm{PSL}(3,3)$ of order $5616$. Any group of the form $\Gamma \cong \pres{H,u}{u^3a^3ua}$ can be decomposed as an amalgamated free product $\Gamma \cong H \ast_{\Z_6} E$ and it follows that $\Gamma$ contains elements of finite order that are not conjugate to any element in $H$. By \cite[(0.4)]{BP}, this in turn implies that if $\mathrm{o}(\beta) = 6$, then the relative presentation $\pres{H,u}{u^3\beta^3 u\beta}$ in \textbf{E2} is not aspherical. Similar remarks apply to the case $\pres{H,u}{u^3\alpha u\alpha^3}$ upon replacing $u$ by $u^{-1}$. Thus only the cases \textbf{E1, E3} remain unresolved for asphericity status. At this point it is not known whether the groups $\pres{a,u}{a^6, u^3a^{\pm 2}ua}$ are finite. We can now state an updated and symmetrized version of the main result from \cite{BBP}.

\begin{Theorem}\emph{\textbf{(\cite[Theorem 4]{BBP},\cite{AE14,GW15})}}\label{Theorem:BBPUpdate} Consider a relative presentation $$\mathcal{Q} = \pres{H,u}{u^3\alpha u\beta },$$ where $1 \neq \alpha ,\beta  \in H$ and suppose that if $\sgp{\alpha ,\beta } \cong \Z_6$, then $\alpha  \neq \beta ^{\pm 2}$ and $\beta  \neq \alpha ^{\pm 2}$. Then $\mathcal{Q}$ is aspherical if and only if none of the following conditions holds:
\begin{enumerate}
\item $\alpha  = \beta ^{\pm 1}$ has finite order;
\item $\alpha  = \beta ^2$ or $\beta  = \alpha ^2$ and $\sgp{\alpha ,\beta } \cong \Z_4$ or $\Z_5$; %$\alpha  = \beta ^2$ and $\mathrm{o}(\beta ) = 4,5$ or $\beta  = \alpha ^2$ and $\mathrm{o}(\alpha ) = 4,5$;
\item $\alpha  = \beta ^3$ or $\beta  = \alpha ^3$ and $\sgp{\alpha ,\beta } \cong \Z_6$;
\item $\{\mathrm{o}(\alpha ), \mathrm{o}(\beta )\} = \{2,3\}$ and $\sgp{\alpha ,\beta } \cong \Z_2 \times \Z_3$;
\item $\frac{1}{\mathrm{o}(\alpha )} + \frac{1}{\mathrm{o}(\beta )} + \frac{1}{\mathrm{o}(\alpha \beta ^{-1})} > 1$ where $1/\infty = 0$.
\end{enumerate}
\end{Theorem}

We now apply Theorem \ref{Theorem:BBPUpdate}, combined with Theorem \ref{Theorem:AsphTransf}, to the relative presentations obtained in Lemma \ref{Lemma:B}, and thence to the cyclic presentations $\mathcal{P}_n(x_0x_jx_kx_l)$. 

\begin{Theorem}\label{Theorem:Asphericity} Assume that $\gcd(n,j,k,l) = 1$. Assume that (B) is true and that both (A) and (C) are false. If $\mathcal{P}_n(x_0x_jx_kx_l)$ is not combinatorially aspherical, then $\mathcal{P}_n(x_0x_jx_kx_l)$ is of type (I) or of type (U). 
\end{Theorem}

\proof By Lemma \ref{Lemma:k=2j} and Theorem \ref{Theorem:DUI}, it suffices to prove the theorem in the case where $k \equiv 2j$ mod $n$ and the parameters $(n,j,k,l)$ satisfy $\gcd(n,j,k,l) = 1$ and both (A) and (C) are false. Let $G = G_n(x_0x_jx_kx_l)$. By Lemma \ref{Lemma:B}, the shift extension $E = G \rtimes_{\theta_G} \Z_n$ admits a relative presentation of the form $\pres{\Z_n,u}{u^3\alpha u \beta}$ where 
$$
\alpha = a^{l-3j}\ \ \mbox{and}\ \ \beta = a^{-l-j}
$$
are such that $\alpha, \beta \neq 1$ and $\alpha \neq \beta^{\pm 1}$ in $\sgp{a} \cong \Z_n$. Theorem \ref{Theorem:AsphTransf} implies that $\mathcal{P}_n(x_0x_jx_kx_l)$ is combinatorially aspherical if and only if the relative presentation $\pres{\Z_n,u}{u^3 \alpha u \beta}$ is aspherical.

It follows from Theorem \ref{Theorem:BBPUpdate} that the relative presentation $\pres{\Z_n,u}{u^3 \alpha u \beta}$ is aspherical unless one of the following conditions is satisfied.
 
 \begin{enumerate}
 \item[(a)] $\alpha = \beta^2$ and $4 \leq \mathrm{o}(\beta) \leq 6$
 \item[(a$'$)] $\beta = \alpha^2$ and $4 \leq \mathrm{o}(\alpha) \leq 6$
 \item[(b)] $\alpha = \beta^{-2}$ and $\mathrm{o}(\beta) = 6$
 \item[(b$'$)] $\beta = \alpha^{-2}$ and $\mathrm{o}(\alpha) = 6$
 \item[(c)] $\alpha = \beta^3$ and $\mathrm{o}(\beta) = 6$
 \item[(c$'$)] $\beta = \alpha^3$ and $\mathrm{o}(\alpha) = 6$
 \item[(d)] $\{\mathrm{o}(\alpha ), \mathrm{o}(\beta )\} = \{2,3\}$ and $\sgp{\alpha ,\beta } \cong \Z_2 \times \Z_3$
 \item[(e)] $\frac{1}{\mathrm{o}(\alpha )} + \frac{1}{\mathrm{o}(\beta )} + \frac{1}{\mathrm{o}(\alpha \beta ^{-1})} > 1$
 \end{enumerate}
Using the fact that $\gcd(n,j,k,l) = 1$, we prove that if any of these conditions hold, then (I) or (U) is satisfied. 

First note that it suffices to consider just one from each of the paired conditions (a)-(a$'$), (b)-(b$'$), and (c)-(c$'$). This is because switching the roles of $\alpha$ and $\beta$ arises from the $\Gamma_n$-transformations
$$
x_0x_jx_kx_l \stackrel{\sigma}{\ra} x_0x_lx_kx_j \stackrel{\tau^2}{\ra} x_kx_jx_0x_l \stackrel{\theta_F^{-k}}{\ra} x_0x_{j-k}x_{-k}x_{l-k} \stackrel{\times-1}{\ra} x_0x_{k-j}x_kx_{k-l} = x_0x_{j'}x_{k'}x_{l'},
$$
where 
$$
\alpha' = a^{l'-3j'} = a^{k-l-3(k-j)} = a^{3j-2k-l} \stackrel{k \equiv 2j}{=} a^{-j-l} = \beta
$$
and 
$$
\beta' = a^{-l'-j'} = a^{l-k-(k-j)} = a^{l-2k+j} \stackrel{k \equiv 2j}{=} a^{l-3j} = \alpha.
$$
The situation (e) does not provide any additional cases because $\sgp{\alpha, \beta}$ is finite cyclic. 

We now consider the cases (a), (b), (c), and (d) individually. All congruences are modulo $n$.

\paragraph{(a) $\alpha = \beta^2$ and $4 \leq \mathrm{o}(\beta) \leq 6$:} Here we have $l - 3j \equiv -2l-2j$ and $4 \leq \mathrm{o}(\beta) = n/\gcd(n,j+l) \leq 6$. That is, we have $$j \equiv 3l\ \ \mbox{and}\ \ n = \mathrm{o}(\beta)\gcd(n,j+l).$$ Since $k \equiv 2j$ we also have $k \equiv 6l$. Now the fact that $1 = \gcd(n,j,k,l) = \gcd(n,3l,6l,l)$ implies that $\gcd(n,l) = 1$ and so $w = x_0x_{3l}x_{6l}x_l$ where $l \in \Z_n^\ast$. Now $j+l \equiv 3l+l$ so $n = \mathrm{o}(\beta)\gcd(n,4l) = \mathrm{o}(\beta)\gcd(n,4)$.

If $\mathrm{o}(\beta) = 4$, then $n = 4\gcd(n,4)$ and it follows that $n=16$, so this is the exemplar for the isolated type (I16). 

If $\mathrm{o}(\beta) = 5$, then $n = 5\gcd(n,4)$ and it follows that $n=5,10$ or $20$. If $n=5$, then $x_0x_3x_6x_1 = x_0x_3x_1x_1$ so we are in the isolated type (I5). If $n=10$, then this is the isolated type (I10). If $n=20$, then this is the isolated type (I20).

If $\mathrm{o}(\beta) = 6$, then $n = 6\gcd(n,4)$ and it follows that $n=24$. This is the unresolved type (U24$'$). 

%\paragraph{\underline{$\beta = \alpha^2$ and $4 \leq \mathrm{o}(\alpha) \leq 6$}:} Here we have $-l-j \equiv 2l-6j$ and $4 \leq \mathrm{o}(\alpha) = n/\gcd(n,l-3j) \leq 6$. That is, we have $$j \equiv 3(2j-l)\ \ \mbox{and}\ \ n = \mathrm{o}(\alpha)\gcd(n,l-3j).$$ Now $k \equiv 2j \equiv 6(2j-l)$ and $l = 2j-(2j-l) \equiv k - (2j-l) \equiv 5(2j-l)$. Thus $1 = \gcd(n,j,k,l) = \gcd(n,3(2j-l),6(2j-l),5(2j-l)) = \gcd(n,2j-l)$ and $w = x_0x_{3u}x_{6u}x_{5u}$ where $u = \gcd(n,2j-l) \in \Z_n^\ast$. Now $l - 3j \equiv 5(2j-l) -9(2j-l)$ so $n = \mathrm{o}(\alpha)\gcd(n,-4(2j-l)) = \mathrm{o}(\alpha)\gcd(n,4)$.
%
%If $\mathrm{o}(\alpha) = 4$, then $n = 4 \gcd(n,4) = 16$ and $x_0x_3x_6x_5$ is in the $\Gamma_{16}$-orbit of $x_0x_1x_2x_7$ so this is the isolated type (I16).
%
%If $\mathrm{o}(\alpha) = 5$, then $n = 5\gcd(n,4)$, so $n=5,10$ or $20$.  If $n=5$ then $w$ is in the $\Gamma_5$ orbit of $x_0x_0x_1x_2$ so this is isolated type (I5). If $n=10$ then $w$ is in the $\Gamma_{10}$ orbit of $x_0x_1x_2x_5$ so this is isolated type (I10). If $n=20$ then $w$ is in the $\Gamma_{20}$ orbit of $x_0x_1x_2x_7$ so this is isolated type (I20).
%
%If $\mathrm{o}(\alpha) = 6$, then $n = 6\gcd(n,4) = 24$ and $x_0x_3x_6x_5$ is in the $\Gamma_{24}$-orbit of $x_0x_3x_6x_1$ so this is the unresolved type (U24$'$).

\paragraph{(b) $\alpha = \beta^{-2}$ and $\mathrm{o}(\beta) = 6$:} In this case $l - 3j \equiv 2j+2l$ and $n = 6\gcd(n,j+l)$. Thus $l \equiv -5j$ and $k \equiv 2j$ so $1 = \gcd(n,j,k,l) = \gcd(n,j,2j,-5j)$ and so $\gcd(n,j) = 1$. Thus $j \in \Z_n^\ast$ and $w = x_0x_jx_{2j}x_{-5j}$. Next, $n = 6\gcd(n,j+l) = 6 \gcd(n,-4j) = 6\gcd(n,4)$ and it follows that $n=24$. So $w = x_0x_jx_{2j}x_{19j}$ is in the unresolved type (U24$''$). 

%\paragraph{\underline{$\beta = \alpha^{-2}$ and $\mathrm{o}(\alpha) = 6$}:} Here we find $l \equiv 7j$ and $k \equiv 2j$ where $j \in \Z_n^\ast$ so $w = x_0x_jx_{2j}x_{7j}$ and $n = 6\gcd(n,l-3j) = 6\gcd(n,4j) = 6\gcd(n,4)$, whence $n=24$ and we have the unresolved type (U24$''$).

\paragraph{(c) $\alpha = \beta^3$ and $\mathrm{o}(\beta) = 6$:} Here we have $l-3j \equiv 3(-l-j)$ so that $4l \equiv 0$, and $n = 6\gcd(n,j+l)$. The fact that $\gcd(n,j,k,l) = 1$ and $k \equiv 2j$ provides that
$$
1 = \gcd(n,j,2j,l) = \gcd(n,j,l) = \gcd(n,j+l,l) = \gcd(l,\gcd(n,j+l)) = \gcd(l,n/6).
$$
Since $4l \equiv 0$ it follows that $\frac{n}{2} \mid 2l$ and so $\frac{n}{6} \mid 2l$. The fact that $\gcd(l,n/6) = 1$ thus implies that $n/6 = 1$ or $2$, so $n = 6$ or $12$.

The element $\alpha = a^{l-3j}$ has order two in $\sgp{a} \cong \Z_n$ so $n = 2\gcd(n,l-3j)$ and so the fact that $6 \mid n$ means that $3 \mid l-3j$ and so $3 \mid l$. %With this, the fact that $\gcd(n,j,k,l) = \gcd(n,j,l) = 1$ implies that $3$ is not a divisor of $j$.

Now suppose that $n = 6$, so that $l \equiv 0$ or $3$. If $l\equiv 0$, then $1 = \gcd(n,j,l) = \gcd(n,j)$ so $j \in \Z_n^\ast$ and $w = x_0x_jx_{2j}x_0$ is a cyclic permutation of $x_0x_0x_jx_{2j}$ as in (I6$''$). If $l \equiv 3$, then $n = 6\gcd(n,j+l) = 6\gcd(n,j+3) = 6$ so $\gcd(n,j+3) = 1$, which implies that $j \equiv 2$ or $4$. Thus $w = x_0x_2x_4x_3 = x_0x_{-4}x_{-2}x_3$ or $x_0x_4x_2x_3$ as in (I6$'$). 

Suppose that $n=12$, so $l \equiv 0,3,6$, or $9$. We have $12 = 6\gcd(n,j+l)$, so $j+l$ is even and the fact that $\gcd(n,j,k,l) = \gcd(n,j,l) = 1$ implies that $j$ and $l$ are both odd, so that $l \equiv 3$ or $9$. Neither $3$ nor $4$ is a divisor of $j+l$, so $j+l \equiv 2$ or $10$. The only arrangements meeting these requirements are $(j,l) \equiv (7,3), (11,3), (1,9)$, or $(5,9)$. In each of these four cases we have $l \equiv 9j$ so $w = x_0x_jx_{2j}x_{9j}$ where $j \in \Z_{12}^\ast$ as in (I12).

%Since $\mathrm{o}(\alpha) = \mathrm{o}(a^{l-3j}) = 2$, we have $2l \equiv 6j \equiv 3k$. Thus:
%\begin{eqnarray*}
%k&\equiv & 3k-2k \equiv 2(l-k),\\
%l &\equiv & 3l-2l \equiv 3(l-k),\ \mbox{and}\\
%j &\equiv & 2j-j \equiv k-j \equiv k-2k \equiv -k \equiv -2(l-k).
%\end{eqnarray*}
%Thus the fact that $\gcd(n,j,k,l) = 1$ implies that $\gcd(n,l-k) = 1$ so $u = l-k \in \Z_n^\ast$. Now the fact that $\mathrm{o}(\beta) = 6$ means that $n = 6\gcd(n,j+l) = 6\gcd(n,l-k) = 6$ and so $w = x_0x_jx_kx_l = x_0x_{4u}x_{2u}x_{3u}$ as in (I6$'$).

%\paragraph{\underline{$\beta = \alpha^3$ and $\mathrm{o}(\alpha) = 6$}:} This analysis is incomplete, but the types (I6$''$) and (I12) both occur here.

%Here we have $-l-j \equiv 3(l-3j)$ and the fact that $\mathrm{o}(\beta) = 2$ implies that $2j+2l \equiv 0$. Thus
%\begin{eqnarray*}
%j&\equiv & 2j-j \equiv k-j,\\
%k &\equiv & 2k-k \equiv j-k \equiv -(k-j),\ \mbox{and}\\
%l &\equiv & 3l-2l \equiv -l-j+9j+2j \equiv -l+2j+7j+(k-j) \equiv 2(k-j) + 7(k-j) \equiv 9(k-j).
%\end{eqnarray*}

\paragraph{(d) $\{\mathrm{o}(\alpha ), \mathrm{o}(\beta )\} = \{2,3\}$ and $\sgp{\alpha ,\beta } \cong \Z_2 \times \Z_3$:} It suffices to consider the case where $\mathrm{o}(\alpha) = \mathrm{o}(a^{l-3j}) = 2$ and $\mathrm{o}(\beta) = \mathrm{o}(a^{-l-j}) = 3$. Here, $2l\equiv 6j$ and $3j \equiv -3j$, so $l \equiv 3l-2l \equiv -3j -6j \equiv -9j$. Coupled with the fact that $k \equiv 2j$, we have $1 = \gcd(n,j,k,l) = \gcd(n,j,2j,-9j)$, so $\gcd(n,j) = 1$ and hence $j \in \Z_n^\ast$ with $w = x_0x_jx_{2j}x_{-9j}$. Now from $\mathrm{o}(\alpha) = \mathrm{o}(a^{l-3j}) = 2$ we have
$$
n = 2 \gcd(n,l-3j) = 2\gcd(n,-12j) = 2\gcd(n,12),
$$
which implies that $8$ is a divisor of $n$ and so with $\mathrm{o}(\beta) = \mathrm{o}(a^{-l-j}) = 3$, it follows that 
$$
n = 3\gcd(n,-l-j) = 3\gcd(n,8j) = 3\gcd(n,8) = 24.
$$
Thus $w = x_0x_jx_{2j}x_{15j}$ as in (I24). This completes the proof.\qed

\bigskip

The following four computational results will help us deal with the isolated types.

\begin{Lemma}\label{Lemma:J4} The group $J_4 = \pres{t,y}{t^4, y^3t^2yt}$ is finite metacyclic of order $272$. The centralizer $\mathrm{Cent}_{J_4}(t)$ is cyclic of order $16$.
\end{Lemma}

\proof Using the GAP commands $\mathrm{Size}, \mathrm{DerivedSubgroup}, \mathrm{AbelianInvariants}$, and $\mathrm{Centralizer}$ we find that $|J_4| = 272$, $J_4/J_4' \cong \Z_{16}$, $J_4' \cong \Z_{17}$, and $\mathrm{Cent}_{J_4}(t) \cong \Z_{16}$. We remark that $J_4$ is the group $J_4(4,1)$ of \cite{BW1}, where it is shown that $J_4$ is a semidirect product $\Z_{17} \rtimes \Z_{16}$. \qed

\begin{Lemma}\label{Lemma:E5} The group $K = \pres{t,u}{t^5, u^3t^2ut}$ is finite metacyclic of order $1100$. The centralizer $\mathrm{Cent}_K(t)$ is an elementary abelian group of order $25$. 
\end{Lemma}

\proof Computing as in Lemma \ref{Lemma:J4}, we find that $|K| = 1100$, $K/K' \cong \Z_{20}$, $K' \cong \Z_{55}$, and $\mathrm{Cent}_{K}(t) \cong \Z_5^2$. \qed 

\begin{Lemma}\label{Lemma:J6} The group $J_6 = \pres{t,u}{t^6, u^3t^3ut^2}$ is finite metacyclic of order $4632$. The centralizer $\mathrm{Cent}_{J_6}(t)$ is cyclic of order $24$.
\end{Lemma}

\proof Here we have $|J_6| = 4632$, $J_6/J_6' \cong \Z_{24}$, $J_6' \cong \Z_{193}$, and $\mathrm{Cent}_{J_6}(t) \cong \Z_{24}$. We remark that $J_6$ is the group $J_6(4,1)$ of \cite{BW1}, where it is shown that $J_6$ is a semidirect product $\Z_{193} \rtimes \Z_{24}$. \qed

\begin{Lemma}\label{Lemma:E6} The group $L = \pres{t,u}{t^6, u^3t^3ut}$ is a nonsolvable group of order $24\,530\,688 = 2^8 \cdot 3^4 \cdot 7 \cdot 13^2$ with second derived subgroup $L''$ isomorphic to the simple group $\mathrm{PSL}(3,3)$ of order $5616$. The centralizer $\mathrm{Cent}_L(t)$ contains an element that does not lie in the normal closure $\<\<t\>\>_L$ of $t$ in $L$. %The centralizer $\mathrm{Cent}_L(t)$ contains an element that is not a power of $t$. 
\end{Lemma}

\proof Although the computations are somewhat lengthy, the order and nonsolvability claims for $L$ can be verified directly in GAP using $\mathrm{Size}, \mathrm{DerivedSubgroup}$, and $\mathrm{GroupId}$ commands, and this was originally done by Williams \cite{GW15}.  The centralizer claim is more difficult and is resistant to brute force application of the $\mathrm{Centralizer}$ command. The group $L$ is the common shift extension for the groups $G_6(x_0x_4x_2x_3)$ of (I6$'$) and $G_6(x_0^2x_1x_2)$ of (I6$''$). One easily checks that 
$$
L \cong E = \pres{a,x}{a^6,x^2axaxa^{-2}} \cong G_6(x_0^2x_1x_2) \rtimes_{\theta} \Z_6
$$
via $a = t^{-1}$ and $x = ua^{-1} = ut$. A direct coset enumeration using $\mathrm{Size}$ returns $|E| = 24\,530\,688$. Now consider the word
\begin{equation}\label{Equation:v}
v = a^3x^3axa^3xa^5xa^5xa^2xa^4xaxa^3xa^2xa^5xa^2x.
\end{equation}
A second coset enumeration reveals that the quotient group $$Q = \pres{a,x}{a^6, x^2axaxa^{-2}, ava^{-1}v^{-1}}$$ has the identical order $|Q| = 24\,530\,688$ and so $v \in \mathrm{Cent}_E(a)$. The fact that $x$ occurs in the word $v$ with exponent sum $14$ implies that the element $v \in E$ is not in the normal closure $\<\<a\>\>_E$ of $a$ in $E$ because $E/\<\<a\>\>_E \cong \pres{x}{x^4} \cong \Z_4$. The element $v \in E$ translates back to $L$ to provide the desired centralizing element. \qed

\paragraph{Remark:} The word $v$ at (\ref{Equation:v}) was discovered using two separate processes in GAP that rely on a stored coset table for the subgroup $\sgp{a} \cong \Z_n$ in $E \cong \pres{a,x}{a^6,x^2axaxa^{-2}}$. This table has columns numbered $1$ to $4\,088\,448$, one for each (right) coset $\sgp{a}w$, with column $1$ corresponding to the coset $\sgp{a}1$. The column corresponding to the coset $\sgp{a}w$ has numbered entries corresponding to the cosets $\sgp{a}wa, \sgp{a}wa^{-1}, \sgp{a}wx, \sgp{a}wx^{-1}$ that arise from the action of the generators $a,x$ of $E$ and their inverses. The first process identifies those columns $j$ for which the entry in the first (and second) row is also $j$, as this corresponds to a coset for which $\sgp{a}v = \sgp{a}va$, so that $v \in \mathrm{Cent}_E(a)$. The second process uses a spanning tree algorithm to extract a word $v$ from the column index $j$. The existing GAP command $$\mathrm{TracedCosetFpGroup(tablename, word, column index)}$$ produces a column index from a word and a column index. (See 47.6-2 in the GAP Manual.) The spanning tree algorithm does the reverse, outputting a (nonempty) word $v$ that solves the equation
$$
j = \mathrm{TracedCosetFpGroup}(\mathrm{tablename}, v, 1)
$$ 
where $j$ is the column index corresponding to a centralizing word found in the first process. The tree search is such that the generators $a$ and $x$ both occur with only positive exponents in the output word $v$ and so that the exponent sum of $x$ in $v$ is as small as possible. The validity of the algorithm is independently verified by the computations in the proof of Lemma \ref{Lemma:E6}.

\begin{Theorem}\label{Theorem:INot} If $\mathcal{P}_n(x_0x_jx_kx_l)$ is of type (I$\,^\ast$), then the shift action by $\Z_n$ on the nonidentity elements of $G = G_n(x_0x_jx_kx_l)$ is not free, so $\mathcal{P}_n(x_0x_jx_kx_l)$ is not combinatorially aspherical. In types (I5), (I6$\'$), and (I6$\"$), $G$ is finite and the shift $\theta_G$ has a nonidentity fixed point.
\end{Theorem}

\proof Let $c = \gcd(n,j,k,l)$ be the primary divisor. If (I$^\ast$) is true, then the parameters $(n/c,j/c,k/c,l/c)$ satisfy the condition (I) and $G$ is the free product $G \cong \ast_{i=1}^c H$ where $H = G_{n/c}(x_0x_{j/c}x_{k/c}x_{l/c})$. Further, as described in \cite[Lemma 5.1]{BShift} for the case of length three relators, the shift on $H$ arises as a restriction of a power of the shift on $G$: $\theta_H = \theta_G^c|_H$. It follows that if $\Z_{n/c}$ does not act freely via the shift on $H$, then $\Z_n$ does not act freely via the shift on $G$. This means that it suffices to prove the theorem under the assumption that the condition (I) is true for the parameters $(n,j,k,l)$. 

%where $\thetaso by Theorem \ref{Theorem:INot}, the cyclic group $\Z_{n/c}$ does not act freely on $H = G_{n/c}(x_0x_{j/c}x_{k/c}x_{l/c})$. Now $G$ is the free product $G \cong \ast_{i=1}^c H$ and $\theta_G^c$ restricts to the shift on $H$: $\theta_G^c|_H = \theta_H$. It follows that if a nonidentity power $\theta_H^p$ fixes a nonidentity element of $H$, then the nonidentity power $\theta_G^{pc}$ fixes that same element, and so $\Z_n$ does not act freely on the nonidentity elements of $G$. 

Taking cases in turn, if $n=5,10, 20$, or $16$, then for type (I$n$) we consider $G = G_n(x_0x_3x_6x_1)$, which has shift extension
$$
E \cong \pres{a,x}{a^n, xa^3xa^3xa^{-5}xa^{-1}} \stackrel{u=xa^3}{\cong} \pres{a,u}{a^n, u^3 a^{-8} u a^{-4}}
$$
and where $G$ is the kernel of a retraction $\nu:E \ra \Z_n$ onto the cyclic group $\sgp{a}$ of order $n$. 

If $n=5$, then $E \cong \pres{a,u}{a^5,u^3a^2ua} \cong K$ as in Lemma \ref{Lemma:E5}, so the centralizer $\mathrm{Cent}_E(a)$ has order $25$. Since $a \in \mathrm{Cent}_E(a)$, the retraction $\nu$ restricts to a surjection of $\mathrm{Cent}_E(a)$ onto $\Z_5$ and so $\mathrm{Fix}(\theta_G) = G \cap \mathrm{Cent}_E(a)$ is cyclic of order $5$.

If $n = 10$ or $20$, then $E \cong \pres{a,u}{a^{n},u^3a^{-8}ua^{-4}} \cong \Z_{n} \ast_{a^{-4} = t} \pres{t,u}{t^5,u^3t^2ut}$ where $\pres{t,u}{t^5,u^3t^2ut} = K$ as in Lemma \ref{Lemma:E5}, so the centralizer $\mathrm{Cent}_K(t)$ has order $25$. Since $t=a^{-4} \in \mathrm{Cent}_K(t)$, the retraction $\nu$ maps $\mathrm{Cent}_K(t)$ onto the subgroup $\sgp{a^4}$ of order $5$ in $\Z_{n}$ and so $\mathrm{Fix}(\theta_G^4) = G \cap \mathrm{Cent}_E(a^4)$ contains a cyclic group of order $5$. 

If $n=16$, then $E \cong \pres{a,u}{a^{16},u^3a^{-8}ua^{-4}} \cong \Z_{16} \ast_{a^{-4} = t} \pres{t,u}{t^4,u^3t^2ut}$ where $\pres{t,u}{t^4,u^3t^2ut} = J_4$ as in Lemma \ref{Lemma:J4}, so the centralizer $\mathrm{Cent}_{J_4}(t)$ has order $16$. Since $t = a^{-4} \in \mathrm{Cent}_{J_4}(t)$, the retraction $\nu$ maps $\mathrm{Cent}_{J_4}(t)$ onto the subgroup $\sgp{a^4}$ of order $4$ in $\Z_{16}$ and so $\mathrm{Fix}(\theta_G^4) = G \cap \mathrm{Cent}_E(a^4)$ contains a cyclic group of order $4$. 

For the types (I6$\'$) and (I6$\"$), consider the group $L = \pres{a,u}{a^6,u^3a^3ua}$ as in Lemma \ref{Lemma:E6}, for which we can choose an element $w \in \mathrm{Cent}_L(a)$ that is not a power of $a$. For $f = 2,5$, there are retractions $\nu^f: L \ra \sgp{a} \cong \Z_6$ satisfying $\nu^f(a) = a$ and $\nu^f(u) = a^f$. Using the rewriting scheme described in \cite[Theorem 2.3]{BShift} together with isomorphisms from Theorem \ref{Theorem:DUI}, we have
\begin{eqnarray*}
\ker \nu^2 & \cong & G_6(u_0u_2u_4u_3) \stackrel{\times -1}{\cong} G_6(u_0u_4u_2u_3)\ \mbox{and}\\
\ker \nu^5 & \cong & G_6(u_0u_5u_4u_0) \stackrel{\tau}{\cong} G_6(u_0^2u_5u_4) \stackrel{\times -1}{\cong} G_6(u_0^2u_1u_2),
\end{eqnarray*}
so these are the groups from (I6$'$) and (I6$''$) respectively. For $f=2,5$, the element $g_f = w\nu^f(w)^{-1}$ is a nontrivial element of $\ker \nu^f \cap \mathrm{Cent}_L(a)$ and so the shifts on the groups $\ker \nu^2 \cong G_6(u_0u_4u_2u_3)$ and $\ker \nu^5 \cong G_6(u_0^2u_1u_2)$ each have nonidentity fixed points.

For type (I12), the group $G = G_{12}(x_0x_1x_2x_9)$ is the kernel of a retraction $\nu: E \ra \Z_{12}$ defined on the shift extension 
$$
E \cong \pres{a,x}{a^{12},xaxaxa^{7}xa^{-9}} \stackrel{u=xa}{\cong} \pres{a,u}{a^{12}, u^3a^6ua^{-10}} \cong \Z_{12} \ast_{a^2=t} \pres{t,u}{t^6, u^3t^3ut}
$$
where $L = \pres{t,u}{t^6, u^3t^3ut}$ is the group from Lemma \ref{Lemma:E6}. An element $w \in \mathrm{Cent}_L(t)$ that is not a power of $t = a^2$ determines a nonidentity element $g = w\nu(w)^{-1} \in G \cap \mathrm{Cent}_E(a^2) = \mathrm{Fix}(\theta_G^2)$.

For type (I24), the group $G = G_{24}(x_0x_1x_2x_{15})$ is the kernel of a retraction $\nu: E \ra \Z_{24}$ defined on the shift extension 
$$
E \cong \pres{a,x}{a^{24},xaxaxa^{13}xa^{-15}} \stackrel{u=xa}{\cong} \pres{a,u}{a^{24}, u^3a^{12}ua^{8}} \cong \Z_{12} \ast_{a^4=t} \pres{t,u}{t^6, u^3t^3ut^2}
$$
where $J_6 = \pres{t,u}{t^6, u^3t^3ut^2}$ is the group from Lemma \ref{Lemma:J6}. Since $t=a^4 \in \mathrm{Cent}_{J_6}(t)$, the retraction $\nu$ maps $\mathrm{Cent}_{J_6}(t)$ onto the subgroup $\sgp{a^4}$ of order $6$ in $\Z_{24}$ and so $\mathrm{Fix}(\theta_G^4) = G \cap \mathrm{Cent}_E(a^4)$ contains a cyclic group of order $4$. \qed

\begin{Theorem}\label{Theorem:FiniteBnotAC} Given $(n,j,k,l)$, let $G = G_n(x_0x_jx_kx_l)$. Assume that (B) is true and (A) and (C) are both false.
\begin{enumerate}
\item[(a)] If $\mathcal{P}_n(x_0x_jx_kx_l)$ is not of type (U$\,^\ast$), then the following are equivalent:
\begin{itemize}
\item[(i)] $\mathcal{P}_n(x_0x_jx_kx_l)$ is combinatorially aspherical;
\item[(ii)] $\Z_n$ acts freely via the shift on the nonidentity elements of $G$;
\item[(iii)] T$\mathcal{P}_n(x_0x_jx_kx_l)$ is not of type (I$\,^\ast$).
\end{itemize}

\item[(b)] The following are equivalent:
\begin{enumerate}
\item[(i)] $G$ is finite;
\item[(ii)] $\mathcal{P}_n(x_0x_jx_kx_l)$ is of type (I) or (U) and the secondary divisor $\gamma = 1$;
\item[(iii)] $\mathcal{P}_n(x_0x_jx_kx_l)$ is of type (I5) or (I6$\'$) or (I6$\"$);
\item[(iv)] The shift $\theta_G$ has a nonidentity fixed point.
\end{enumerate}
%\item[(b)] If (I5) is true then $G$ is a non-nilpotent metacyclic group of order $220$.
%\item[(c)] The groups in types (I6$\'$) and (I6$\"$) are non-isomorphic non-solvable groups of order $4\,088\,448$.
\end{enumerate}
\end{Theorem}

\proof For part (a), assume that $\mathcal{P}_n(x_0x_jx_kx_l)$ is not of type (U$^\ast$). The implications (i) $\Rightarrow$ (ii) $\Rightarrow$ (iii) follow from Theorems \ref{Theorem:Free} and \ref{Theorem:INot} respectively. If $\mathcal{P}_n(x_0x_jx_kx_l)$ is not of type (I$^\ast$) or (U$^\ast$), then $\mathcal{P}_{n/c}(x_0x_{j/c}x_{k/c}x_{l/c})$ is combinatorially aspherical by Theorem \ref{Theorem:Asphericity}, and so $\mathcal{P}_n(x_0x_jx_kx_l)$ is combinatorially aspherical by Lemma \ref{Lemma:FreeProduct}. This shows that (iii) $\Rightarrow$ (i) and so completes the proof of part (a).

For part (b), assume first that $G$ is finite. The primary and secondary divisors $c = \gcd(n,j,k,l)$ and $\gamma$ are both equal to one by Corollary \ref{Corollary:Primary} and Theorem \ref{Theorem:secondary}, respectively. Since (A) is false, it follows that $n \neq 1$ and so $\mathcal{P}_n(x_0x_jx_kx_l)$ is not combinatorially aspherical by Lemma \ref{Lemma:FiniteAsph}. That (I) or (U) must be true then follows from Theorem \ref{Theorem:Asphericity}. This shows that (i) $\Rightarrow$ (ii). The secondary divisor calculations in Example \ref{Example:IsolatedExamples} show that (ii) $\Rightarrow$ (iii). That (iii) $\Rightarrow$ (iv) and (iii) $\Rightarrow$ (i) follow directly from Theorem \ref{Theorem:INot}. And finally, if the shift $\theta_G$ has a nonidentity fixed point, then the primary and secondary divisors are both equal to one by Theorem \ref{Theorem:secondary} and $\mathcal{P}_n(x_0x_jx_kx_l)$ is not combinatorially aspherical by Theorem \ref{Theorem:Free}, whence (I) or (U) is true by Theorem \ref{Theorem:Asphericity}. Thus (iv) $\Rightarrow$ (ii) and the proof is complete. \qed

\section{Concluding Remarks}\label{Section:Conclusion}

In this concluding section we distill the proofs of the main theorems \ref{Theorem:AsphericityClassification}-\ref{Theorem:Finite=Fix} from the results in Sections \ref{Section:SmallCancellation}-\ref{Section:FTF}. This includes a discussion of the finite groups $G_n(x_0x_jx_kx_l)$ that occur in Theorem \ref{Theorem:FinitenessClassification}. We conclude with a discussion of the unresolved cases from the asphericity classification Theorem \ref{Theorem:AsphericityClassification}.

Avoiding the unresolved situation (U$^\ast$), the proof of Theorem \ref{Theorem:AsphericityClassification} is a consequence of Theorems \ref{Lemma:SmCanc}, \ref{Lemma:Decomposable}, \ref{Theorem:MZ}(a), \ref{Theorem:Asphericity}, and \ref{Theorem:INot}. If $k \not \equiv 0$ or $j \not \equiv l$ mod $n$, then the presentation $\mathcal{P} = \mathcal{P}_n(x_0x_jx_kx_l)$ has no proper power relators, in which case combinatorial asphericity implies that after removing freely redundant relators from the two-dimensional model of $\mathcal{P}$, the group $G_n(x_0x_jx_kx_l)$ is the fundamental group of an aspherical two-complex, and so is torsion-free with geometric dimension at most two.

The finiteness classification Theorem \ref{Theorem:FinitenessClassification} follows from Lemma \ref{Lemma:FiniteAsph}(b) and Theorems \ref{Theorem:FiniteCnotA}, \ref{Theorem:MZ}(b), and \ref{Theorem:FiniteBnotAC}(b). These results also detail the structural claims for the finite groups that occur, as for example in Theorem \ref{Theorem:MZ}(c), with the exception of those occurring in types (I5), (I6$'$), and (I6$''$). The group $G_5(x_0x_3x_1^2)$ of (I5) is a normal subgroup of index $5$ in its shift extension, which is the group $K$ of Lemma \ref{Lemma:E5}. Thus $G_5(x_0x_3x_1^2)$ is metacyclic of order $220$ and in fact is a semidirect product $\Z_5 \rtimes \Z_{44}$. As a nonabelian group with cyclic abelianization, this group is not nilpotent. As in the proof of Lemma \ref{Lemma:E6}, the groups $$G_1 = G_6(x_0x_4x_2x_3)\ \mbox{and}\ G_2 = G_6(x_0^2x_1x_2)$$ of (I6$'$) and (I6$''$) share the common shift extension $L = \pres{a,u}{a^6, u^3a^3ua}$ of order $24\,530\,688$, which has second derived subgroup $L''$ isomorphic to the simple group $\mathrm{PSL}(3,3)$ of order $5616$. Both $G_1$ and $G_2$ therefore are nonsolvable of order $4\,088\,448$ and contain $L'' \cong \mathrm{PSL(3,3)}$. Calculations in GAP show that the abelianizations $G_1^\mathrm{ab} \cong \Z_8$ and $G_2^\mathrm{ab} \cong \Z_7 \oplus \Z_8$ are not isomorphic, so $G_1 \not \cong G_2$. 

For the proof of Theorem \ref{Theorem:AsphericalFree}, since cyclic presentations with positive relators are orientable, Theorem \ref{Theorem:Free} provides that combinatorial asphericity always implies free action. For the converse, and avoiding the unresolved type (U$^\ast$), in each nonaspherical case we have shown that the shift action is not free. See Theorems \ref{Lemma:Decomposable}, \ref{Theorem:MZ}(a), and \ref{Theorem:INot}.

For the proof of Theorem \ref{Theorem:FiniteFix}, for each finite group $G_n(x_0x_jx_kx_l)$, either the shift action is trivial, as when $n=1$ and in Lemma \ref{Lemma:FiniteAsph}(b) and Theorem \ref{Theorem:FiniteCnotA}, or else we have shown that the shift has a nonidentity fixed point, as in Theorems \ref{Theorem:MZ}(b) and \ref{Theorem:INot}. If (C) is true and $\gcd(n,2k) = 1$, then (A) is false and $\mathrm{Fix}(\theta_G) \neq 1$ by Theorem \ref{Theorem:CFix}. Conversely, if $\mathrm{Fix}(\theta_G) \neq 1$, then either $n=1$ so that $G \cong \Z_4$ is finite (Lemma \ref{Lemma:FiniteAsph}(b)), (C) is true and $\gcd(n,2k) = 1$ (Theorem \ref{Theorem:CFix}), or $G$ is finite (Theorems \ref{Theorem:MZ}(b), \ref{Theorem:FiniteBnotAC}(b)).

In order to prove Theorem \ref{Theorem:Finite=Fix}, suppose that $\mathrm{Fix}(\theta_G) \neq 1$ where $G$ is infinite (for otherwise we may take $G = H$). By Theorem \ref{Theorem:FiniteFix}, (C) is true and $\gcd(n,2k) = 1$. By Lemma \ref{Lemma:Cdecomp}, there exists an integer $p$ such that the shift extension $E = G \rtimes_{\theta_G} \Z_n$ admits a presentation of the form 
$$
E \cong %\pres{a,z}{a^n, z^2a^{k-2p}z^2a^{-k-2p}} 
\pres{a,u,z}{a^n, ua^kua^{-k}, ua^{2p} = z^2}
$$
where $G = \ker \nu$ for a retraction $\nu:E \ra \Z_n = \sgp{a}$ satisfying $\nu(a) = a$, $\nu(u) = 1$, and $\nu(z) = a^p$. And finally, since $G$ is infinite, Theorem \ref{Theorem:FiniteCnotA} implies that $\gcd(n,p) \neq 1$.  Now $E \cong D \ast_{\sgp{
ua^{2p} = z^2}} \sgp{z}$ where $D = \pres{a,u}{a^n, ua^kua^{-k}} \cong \Z_n \times \Z_2$. As in \cite[Theorem 2.3]{BShift}, it follows that $G = \ker \nu$ contains the finite group $H = \ker \nu|_D \cong G_n(u_0u_k) \cong \Z_2$ and so $\mathrm{Fix}(\theta_G)$ contains $\mathrm{Fix}(\theta_H) = \Z_2 \neq 1$. It remains to show that $\mathrm{Fix}(\theta_G) = \mathrm{Fix}(\theta_H)$. So assume that $g \in \mathrm{Fix}(\theta_G) = \mathrm{Cent}_E(a) \cap G$. By the Centralizer Lemma \ref{lemma:cent}, either $g \in D \cap G = H$ or else there exists $d \in D$ such that $a \in d\sgp{ua^{2p}}d^{-1}$. In the former case we have $g \in \mathrm{Fix}(\theta_H)$, as desired, and the latter case cannot occur due to the fact that $\gcd(n,p) \neq 1$, so $\nu(ua^{2p})$ generates a proper subgroup of $\nu(E) = \sgp{\nu(a)} = \Z_n$. Thus $\mathrm{Fix}(\theta_G) = \mathrm{Fix}(\theta_H)$ and the proof is complete.

We conclude with a discussion of the unresolved type (U$^\ast$), which boils down to consideration of the two exemplar presentations from Table \ref{Table:U}: 
$$
\mathcal{P}_{24}(x_0x_3x_6x_1)\ \mbox{and}\ \mathcal{P}_{24}(x_0x_1x_2x_{19}).
$$
In both cases, the secondary divisor is equal to four so by Theorem \ref{Theorem:secondary}, the groups defined by these presentations are infinite and their shifts are fixed point free. However it is unknown whether these presentations are combinatorially aspherical or whether either shift action is free. It would be of interest to determine, for example, whether these groups are torsion-free. At this point we can note that the groups $G_1 = G_{24}(x_0x_3x_6x_1)$ and $G_2 = G_{24}(x_0x_1x_2x_{19})$ are not isomorphic. Routine calculations in GAP provide that they abelianize to $G_1/G_1' \cong \Z^3 \oplus \Z_5^2 \oplus \Z_7$ and $G_2/G_2' \cong \Z^3 \oplus \Z_3 \oplus \Z_{73}$. It seems likely that $G_1$ and $G_2$ are commensurable, although their shift extensions 
$$
E_1 \cong \pres{a,x}{a^{24}, xa^3xa^3xa^{-5}xa^{-1}}\ \mbox{and}\ E_2 \cong \pres{a,x}{a^{24}, xaxaxa^{17}xa^{-19}}
$$
are not isomorphic, having second derived quotients $E_1'/E_1'' \cong \Z^9 \oplus \Z_5^8 \oplus \Z_7^4$ and $E_2'/E_2'' \cong \Z^9 \oplus \Z_3^4 \oplus \Z_{73}^4$. %Let $\mathcal{P} = \mathcal{P}_n(x_0x_jx_kx_l)$ and $G = G_n(x_0x_jx_kx_l)$. 

%
%
%
%

%\begin{thebibliography}{99}
%
%\bibitem{NeumannRoot} B.~H.~Neumann, Adjunction of elements to groups, J.~London Math.~Soc.~{\bf 18}, (1943), 4-11.
%\end{thebibliography}
%
\end{document}